\documentclass[11pt]{article}

\usepackage[margin=1in]{geometry}
\usepackage{amsmath,amssymb,amsthm}
\usepackage{graphicx}
\usepackage{booktabs}
\usepackage{array}
\usepackage[dvipsnames,xcdraw]{xcolor}
\usepackage{hyperref}
\usepackage{placeins}
\usepackage{tikz-cd}

\numberwithin{equation}{section}

\theoremstyle{plain}
\newtheorem{theorem}{Theorem}[section]
\newtheorem{conjecture}[theorem]{Conjecture}

\newtheorem{proposition}[theorem]{Proposition}
\newtheorem{lemma}[theorem]{Lemma}
\newtheorem{corollary}[theorem]{Corollary}
\newtheorem{hypothesis}[theorem]{Hypothesis}

\theoremstyle{definition}

\theoremstyle{definition}
\newtheorem{remark}{Remark}[section]

\newcommand{\ApmMACaption}{Raw values are shown in gray every $500$ in $x$, and the colored curve is the symmetric moving average over $[x-2500,x+2500]$.}

\newcommand{\Dset}{\mathcal D}

\newcommand{\F}{\mathbb F}

\newcommand{\Leg}[2]{\left(\frac{#1}{#2}\right)}

\newcommand{\Z}{\mathbb Z}

\newenvironment{red}{\relax\color{red}}{\relax}
\newenvironment{blue}{\relax\color{blue}}{\hspace*{.5ex}\relax}

\newcommand{\ber}{\begin{red}}
\newcommand{\er}{\end{red}}
\newcommand{\beb}{\begin{blue}}
\newcommand{\eb}{\end{blue}}
\newcommand{\Scalar}{\mathrm{sc}}
\newcommand{\Split}{\mathrm{sp}}
\newcommand{\Ramified}{\mathrm{ram}}

\DeclareMathOperator{\Ind}{Ind}
\DeclareMathOperator{\Li}{Li}
\DeclareMathOperator{\PSL}{PSL}
\DeclareMathOperator{\SL}{SL}
\DeclareMathOperator{\Tr}{Tr}

\usepackage{authblk}

\title{
Character-Twisted Sums of Real Quadratic Class Numbers \\ 
and a Selberg-Type Conjecture
}

\author{Chul-hee Lee, Kyu-Hwan Lee, and Jinsung Park}

\date{\today}

\begin{document}

\maketitle

\begin{abstract}
Extending the arithmetic Selberg-zeta function framework studied by Sarnak \cite{Sarnak1982} and Hashimoto \cite{Hashimoto2007}, we study Dirichlet-character weighted sums of narrow class numbers of real quadratic orders and derive explicit formulas according to the scalar, split, ramified, and inert reductions of primitive hyperbolic classes.
Our first main result establishes asymptotic formulas for the class number sums associated with the three noninert reduction types and shows that, after multiplication by the corresponding nontrivial even-character weights, their leading terms cancel exactly for every fixed prime $p\ge 5$.
Computational experiments across the sampled ranges of levels and geodesic lengths support a substantially stronger bound  for the sum over nontrivial even-characters and exhibit systematic cancellation among the scalar, split, and ramified contributions.
If this conjectural bound persists over arbitrarily long polynomial ranges, the resulting exceptional eigenvalues-free regions approach the Selberg threshold $1/4$, providing a possible arithmetic route toward a Selberg-type spectral gap.
\end{abstract}



\section{Introduction}
\label{sec:introduction}

The Selberg zeta function provides a powerful bridge between the geometry of closed geodesics and the spectral theory of automorphic forms, with additional arithmetic significance when the underlying surface is arithmetic.
For arithmetic surfaces such as the modular surface, this connection becomes especially explicit, since primitive closed geodesics can themselves be described in terms of classical arithmetic data. We begin by recalling this  arithmetic description for the modular surface.

Let
\[
\Gamma(1):=\SL_2(\mathbb Z),
\qquad
Y(1):=\Gamma(1)\backslash\mathbb H,
\]
where $\mathbb H$ is the upper half-plane.  Let
\[
\Dset
=
\{d\in\mathbb Z_{>0}:d\equiv0,1 \!\!\! \pmod4,\ d\text{ is not a square}\}.
\]
For $d\in \Dset$, let \(\mathcal O_d\) be the real quadratic order of discriminant \(d\), and let \(\epsilon_d\) denote its fundamental totally positive unit. 
We write
\begin{equation}\label{eq:Pell-data}
\epsilon_d
=
\frac{t_d+u_d\sqrt d}{2}>1,
\qquad
t_d^2-du_d^2=4,
\end{equation}
and denote by $h(d)$ the narrow class number of $\mathcal O_d$.
The classical correspondence between indefinite binary quadratic forms and primitive hyperbolic conjugacy classes gives, for each discriminant $d$, a collection of primitive closed geodesics on the modular surface $Y(1)$ with common length and norm
\[
\ell(d)=2\log\epsilon_d,
\qquad
N_d=\epsilon_d^2.
\]
We refer to Section~1 of \cite{Sarnak1982} for details of this correspondence. 
Consequently, the Selberg zeta function of \(Y(1)\) can be expressed in terms of
the arithmetic data \(\epsilon_d\) and \(h(d)\); see \eqref{eq:level-one-class-number-log-derivative}. For the corresponding
character-twisted formula for \(\Gamma_0(p)\), whose definition we recall below, see \eqref{eq:VZ-character-zeta} and \eqref{eq:zeta-log-derivative-class-number}.

Exploiting this correspondence between arithmetic and geometry, \cite{Sarnak1982} and \cite{Hashimoto2007} studied several problems concerning the distribution of narrow class numbers of real quadratic fields.
The starting point of the present paper is to extend this framework by introducing nontrivial Dirichlet-character weights.
This naturally brings the congruence subgroup $\Gamma_0(p)$ into consideration, since Dirichlet characters arise through its Borel quotient modulo an odd prime $p$.
{\em A principal motivation for introducing these character weights is to detect cancellations that may reveal hidden symmetries in the distribution of narrow class numbers.} Throughout the remainder of this paper, “class number” refers to the narrow class number.

Let us recall that the standard principal congruence subgroup is defined by
\[
\Gamma(p)
:=
\{\gamma\in\Gamma(1):\gamma\equiv I\pmod p\}.
\]
Reduction modulo $p$ gives
\[
\widetilde G_p
:=
\Gamma(1)/\Gamma(p)
\simeq
\SL_2(\mathbb F_p).
\]
Let
\[
\widetilde B_p:=\Gamma_0(p)/\Gamma(p),
\]
where
\[
\Gamma_0(p)
:=
\left\{
\begin{pmatrix}
a & b\\
c & d
\end{pmatrix}
\in \Gamma(1)
:
c\equiv 0 \pmod p
\right\},
\]
and let $\widetilde U_p$ be its unipotent radical.  Then
\[
\widetilde B_p/\widetilde U_p
\simeq
\mathbb F_p^\times.
\]
Thus every Dirichlet character
\[
\chi:\mathbb F_p^\times\longrightarrow\mathbb C^\times
\]
induces a character $\widetilde\chi$ of $\widetilde B_p$. 
We denote the corresponding induced representation of $\widetilde{G}_p$ by
\begin{equation}\label{eq:induced-rep}
\widetilde\rho_\chi
:=
\Ind_{\widetilde B_p}^{\widetilde G_p}\widetilde\chi.
\end{equation}

For an even Dirichlet character, i.e., $\chi(-1)=1$, the representation
$\widetilde\rho_\chi$ is trivial on the center $\{\pm I\}$ and
therefore descends to
\[
G_p
:=
\widetilde G_p/\{\pm I\}
\simeq
\PSL_2(\mathbb F_p).
\]
Let $B_p$ and $U_p$ denote the images of $\widetilde B_p$ and
$\widetilde U_p$ in $G_p$, and put
\[
H_p
:=
B_p/U_p
\simeq
\mathbb F_p^\times/\{\pm1\}.
\]
Thus the even Dirichlet characters modulo $p$ are precisely the
characters of $H_p$, and for such $\chi$ we write
\[
\rho_\chi
:=
\Ind_{B_p}^{G_p}\chi.
\]

For $d \in \Dset$, let $\gamma\in \Gamma(1)$ be a primitive hyperbolic element corresponding to $d$; see \eqref{eq:Q-matrix}.
The characteristic polynomial is $x^2 - t_d x + 1$, and  its discriminant is
\[
  \operatorname{disc}(x^2 - t_d x + 1) = t_d^2 - 4 = d\,u_d^2.
\]
The Pell equation \eqref{eq:Pell-data}
yields a fourfold classification of the reduction
\(\widetilde{\gamma}\in\SL_2(\mathbb F_p)\) of the primitive
hyperbolic class.  We call it \emph{scalar}, \emph{split},
\emph{ramified}, or \emph{inert} according as
\begin{equation}
\begin{array}{c|c}
\text{type} & \text{Pell condition}\\ \hline
\text{scalar}    & p\mid u_d,\\
\text{split}     & p\nmid u_d,\ \left(\frac{d}{p}\right)=1,\\
\text{ramified}  & p\nmid u_d,\ p\mid d,\\
\text{inert}     & p\nmid u_d,\ \left(\frac{d}{p}\right)=-1.
\end{array}
\label{eq:types}
\end{equation}
Hereafter, we use $\Scalar$, $\Split$, and $\Ramified$ to denote the scalar, split, and ramified types, respectively, in mathematical expressions.
Equivalently, the four types are distinguished by the splitting behavior
of the characteristic polynomial of \(\widetilde\gamma\). Section~\ref{sec:class-number-geodesics}
interprets them through the action on
\(\mathbb P^1(\mathbb F_p)\) and derives the corresponding
character-weighted fixed-line formulas, including the summed identities
\eqref{eq:even-local-type-weight} and \eqref{eq:odd-local-type-weight}.

For a Dirichlet character $\chi$ modulo $p$, let $A_{p,\chi}(x)$ be
the class-number sum defined in \eqref{eq:def-A-p-chi}.  We write
\[
A_p^+(x)
:=
\sum_{\substack{\chi\bmod p\\
\chi(-1)=1,\ \chi\neq \mathbf 1}}
A_{p,\chi}(x),
\qquad
A_p^-(x)
:=
\sum_{\substack{\chi\bmod p\\
\chi(-1)=-1}}
A_{p,\chi}(x),
\]
and
\[
A_p(x):=A_p^+(x)+A_p^-(x).
\]
The formulas of Section~\ref{sec:class-number-geodesics} make the even character sum accessible to computation:
the sum $A_p^+(x)$ over all nontrivial even characters is expressed explicitly according to whether the reduction type is scalar, split, ramified, or inert, with no contribution from the inert case, whose weight is zero.
Thus
$A_p^+(x)$ can be computed directly from the Pell and class-number data.  The
experimental computations in Section~\ref{sec:experimental-setup} reveal substantial cancellation among the scalar, split, and ramified contributions and consistently suggest the following conjecture.

\begin{conjecture}\label{conj:ap-experimental}
For every \(\epsilon>0\), 
\[
   A_p^+(x)\ll_{\epsilon} (p x)^{1+\epsilon},
\]
uniformly for odd primes \(p\) and \(x\geq 2\).
\end{conjecture}

Although Conjecture \ref{conj:ap-experimental} concerns a uniform bound of size \(px\), the numerical plots below use the normalization \(x\sqrt p\). This provides a convenient reference scale for comparing the observed fluctuations across different primes and displaying the cancellation between reduction-type contributions. It should not be interpreted as a uniform upper bound at every sharp cutoff. As explained in Remark \ref{re:Arithmetic jumps and the conjectural bound}, ramified and scalar reductions can produce larger arithmetic jumps, while coarse sampling and averaging may conceal the resulting jumps. We therefore retain the \(x\sqrt p\)-normalization for the numerical illustrations while formulating the uniform conjecture on the \(px\)-scale.

Our main result in Section~\ref{sec:cancellation-thm}  gives separate class number asymptotics for
the three contributing reduction types. For every fixed prime \(p\ge5\), as
\(x\to\infty\),
\[
\begin{aligned}
 \sum_{d\in \Dset_{0,p}^{\Split}(x)} h(d)
 &\sim \frac{p-3}{2(p-1)}\operatorname{Li}(x^2),\\
 \sum_{d\in \Dset_{0,p}^{\Ramified}(x)} h(d)
 &\sim \frac{2}{p}\operatorname{Li}(x^2),\\
 \sum_{d\in \Dset_{0,p}^{\Scalar}(x)} h(d)
 &\sim \frac{2}{p(p^2-1)}\operatorname{Li}(x^2),
\end{aligned}
\]
where \(\operatorname{Li}(x)=\displaystyle\int_2^x \frac{dt}{\log t}\) denotes the logarithmic integral.
Here the superscript \(T\) in $ \Dset^T_{0,p}(x)$ denotes restriction to discriminants of
reduction type \(T\), as defined in \eqref{eq:types}, with the condition
\(\epsilon_d\leq x\) imposed. The notation is defined formally
in Section~\ref{sec:experimental-setup}.
These formulas follow from applying the Chebotarev-type prime geodesic theorem to the reduction-type decomposition in Section~\ref{sec:class-number-geodesics}.
After the corresponding nontrivial even-character weights in \eqref{eq:even-local-type-weight} are inserted,
the leading coefficients satisfy
\[
 -\frac{p-3}{p-1}
 +\frac{p-3}{p}
 +\frac{p-3}{p(p-1)}=0,
\]
and hence
\[
 A_p^+(x)=o\!\left(\operatorname{Li}(x^2)\right).
\]
This cancellation also has a natural representation-theoretic explanation. The representation
\begin{equation}\label{eq:def-rho-Dir}
\rho_{\mathrm{Dir}}
=
\bigoplus_{\substack{\chi\bmod p\\ \chi(-1)=1,\ \chi\ne\mathbf 1}}
\operatorname{Ind}_{B_p}^{G_p}\chi
\end{equation}
contains no trivial \(G_p\)-constituent, since Frobenius reciprocity gives
\[
\operatorname{Hom}_{G_p}
\left(\mathbf 1,\operatorname{Ind}_{B_p}^{G_p}\chi\right)
\simeq
\operatorname{Hom}_{B_p}(\mathbf 1,\chi)
=0
\]
for every nontrivial \(\chi\). Its character therefore has zero average over \(G_p\), so the corresponding Chebotarev main term vanishes. See Theorem \ref{thm:fixed-level-cancellation}
for further details.

Section~\ref{sec:local-type-contributions} begins
by examining the individual contributions from each reduction type,
whose leading terms are determined in
Section~\ref{sec:cancellation-thm}.
Figure~\ref{fig:split-ramified} illustrates the striking
symmetry between the split and ramified contributions
and the resulting substantial cancellation.
We then subtract the main terms and examine the remaining
errors after normalization by \(x\sqrt p\).
Figures~\ref{fig:survivor-split}--\ref{fig:survivor-scalar-sigma}
display these fluctuations, revealing differences among
the reduction types and the effects of the scalar contributions
within the sampled range.

Section~\ref{sec:dirichlet-spectral-gap} establishes a quantitative connection between the arithmetic bound in Conjecture~\ref{conj:ap-experimental} and the exceptional spectrum of the following automorphic space. For a nontrivial even Dirichlet character \(\chi\) modulo \(p\), let
\[
%\label{eq:def-Dirichlet-sector}
\mathcal H_\chi
:=
L^2_{\mathrm{cusp}}
\bigl(\Gamma_0(p)\backslash\mathbb H,\chi\bigr)
\]
denote the Hilbert space of square-integrable cuspidal functions
$f:\mathbb H\to\mathbb C$ satisfying
\[
f(\gamma z)
=
\chi(a\bmod p)\,f(z),
\qquad
\gamma=
\begin{pmatrix}
a&b\\ c&d
\end{pmatrix}
\in\Gamma_0(p).
\]
We then define the nontrivial even Dirichlet sector by
\[
\mathcal H_{\mathrm{Dir}}(p)
:=
\bigoplus_{\substack{\chi\bmod p\\
                     \chi(-1)=1,\ \chi\neq \mathbf 1}}
\mathcal H_\chi.
\]
For $0<r<1/2$, define
\[
  N_{\mathrm{Dir}}(r;p)
  :=
  \sum_{\substack{\chi\bmod p\\
  \chi(-1)=1,\ \chi\ne \mathbf 1}}
  \#\left\{
  j:
  0<\lambda_{\chi,j}
  \le \frac14-r^2
  \right\},
%\label{eq:densityN}
\]
where the $\lambda_{\chi,j}$ are the cuspidal 
eigenvalues of the Laplacian acting on $\mathcal H_\chi$, counted with multiplicity.
Conjecture~\ref{conj:ap-experimental}, together with
Proposition~\ref{prop:conditional-dirichlet-gap},
motivates the following spectral conjecture.

\begin{conjecture}\label{conj:Dirichlet-Selberg}
For every $0<r<1/2$, there exists $p_0=p_0(r)$ such that
\[
N_{\mathrm{Dir}}(r;p)=0
\]
for every odd prime $p\ge p_0$.
\end{conjecture}

The classical Selberg eigenvalue conjecture asserts that every nonzero
cuspidal eigenvalue of the Laplacian for a congruence subgroup of
$\SL_2(\mathbb Z)$ is at least $1/4$; see Selberg~\cite{Selberg1965}.
Conjecture~\ref{conj:Dirichlet-Selberg} is a weak form of the
corresponding assertion restricted to the nontrivial even Dirichlet
sector. {\em Our work connects this spectral problem to the character-weighted class-number sums, opening a possible new direction for investigating exceptional eigenvalues.}

\section{Class numbers, primitive geodesics, and character weights}
\label{sec:class-number-geodesics}

We continue to use the notation and conventions introduced in
Section~\ref{sec:introduction}.  
We shall use both reduction maps, related by the commutative diagram
\begin{equation}
\begin{tikzcd}[column sep=large, row sep=large]
\Gamma(1)
  \arrow[r, "\widetilde{\pi}_p"]
  \arrow[dr, swap, "\pi_p"]
&
\widetilde{G}_p
  \arrow[d, "q_p"]
\\
& G_p ,
\end{tikzcd}
\end{equation}
where \(q_p\colon \widetilde{G}_p\to G_p\) is the quotient map. Thus
\[
\pi_p=q_p\circ\widetilde{\pi}_p.
\]
For $\gamma\in\Gamma(1)$, write
\[
\widetilde\gamma:=\widetilde\pi_p(\gamma),
\qquad
\bar\gamma:=\pi_p(\gamma).
\]
Thus $\widetilde\gamma$ retains the central sign in
$\SL_2(\mathbb F_p)$, whereas $\bar\gamma$ denotes its projective
image in $\PSL_2(\mathbb F_p)$.

\subsection{Selberg zeta functions associated to Dirichlet characters}
\label{subsec:selberg-zeta-class-number}

Let \(\Gamma=\Gamma(1)\) or \(\Gamma_0(p)\), and let \(\rho:\Gamma\to U(V)\) be a finite-dimensional unitary representation. The Selberg zeta function twisted by \(\rho\) is defined by
\[
Z_\Gamma(s,\rho)
=
\prod_{\gamma\in \mathcal P(\Gamma)}
\prod_{n=0}^{\infty}
\det\!\left(I_V-\rho(\gamma)N(\gamma)^{-s-n}\right),
\qquad \Re s>1.
%\tag{2.2}
\]
Here \(\mathcal P(\Gamma)\) denotes the set of primitive hyperbolic conjugacy classes of \(\Gamma/\{\pm I\}\), represented by their positive-trace lifts to \(\Gamma\), and \(N(\gamma)=e^{\ell_\gamma}\), where \(\ell_\gamma\) is the corresponding geodesic length. The product converges absolutely for \(\Re s>1\). For the trivial one-dimensional representation \(\mathbf 1\), we write \(Z_\Gamma(s):=Z_\Gamma(s,\mathbf 1)\). In particular, the untwisted Selberg zeta function for the modular surface \(Y(1)\) has the arithmetic expression
\begin{equation}
\frac{Z'_{\Gamma(1)}(s)}
     {Z_{\Gamma(1)}(s)}
=
\sum_{d\in \Dset}
\sum_{j\geq1}
h(d) \,
\frac{2\log\epsilon_d}
     {1-\epsilon_d^{-2j}}
\, \epsilon_d^{-2js},
\qquad
\Re s>1.
\label{eq:level-one-class-number-log-derivative}
\end{equation}
We refer to 
\cite[(1.2)--(1.3)]{Hashimoto2007} for this.
The factor \(h(d)\) in \eqref{eq:level-one-class-number-log-derivative} arises from grouping primitive hyperbolic conjugacy classes by discriminant. Indeed, the classical correspondence identifies the narrow ideal classes of \(\mathcal O_d\) with the proper equivalence classes of primitive indefinite binary quadratic forms of discriminant \(d\), and these in turn correspond to primitive hyperbolic conjugacy classes. There are therefore exactly \(h(d)\) such classes for each \(d\). Since they all have expanding eigenvalue \(\epsilon_d\), they contribute the same summand to the logarithmic derivative of the Selberg zeta function. Summing their contributions produces the factor \(h(d)\) in \eqref{eq:level-one-class-number-log-derivative}.

The character $\chi$ of $\mathbb F_p^\times$ defines a character of
$\Gamma_0(p)$ by
\[
\chi_{\Gamma_0(p)}
\begin{pmatrix}
a&b\\ c&d
\end{pmatrix}
:=
\chi(a\bmod p).
\]
Since $ad\equiv1\pmod p$, it is related to the usual nebentypus
convention by
\[
\chi(a)=\chi^{-1}(d).
\]
Thus, in the standard convention based on the lower-right entry,
$\chi_{\Gamma_0(p)}$ has nebentypus $\chi^{-1}$.
Equivalently, it is obtained by composing
\[
\Gamma_0(p)
\longrightarrow
\widetilde B_p
\longrightarrow
\widetilde B_p/\widetilde U_p
\simeq\mathbb F_p^\times
\]
with $\chi$.
Recall that
$\widetilde\rho_\chi
=\operatorname{Ind}_{\widetilde B_p}^{\widetilde G_p}
\widetilde\chi$
is a representation of $\widetilde G_p$. Its pullback to
$\Gamma(1)$ by $\widetilde\pi_p$ satisfies
\[
 \widetilde\rho_\chi\circ\widetilde\pi_p
 \simeq
 \operatorname{Ind}_{\Gamma_0(p)}^{\Gamma(1)}
 \chi_{\Gamma_0(p)}.
\]

For even $\chi$, these Euler products are the ordinary weight-zero
Selberg zeta functions, and the Venkov--Zograf induction identity
in \cite{VenkovZograf}
gives
\begin{equation}\label{eq:VZ-character-zeta}
 \begin{aligned}
 Z_{\Gamma_0(p)}(s,\chi_{\Gamma_0(p)})
 &=
 Z_{\Gamma(1)}
 \left(s,\operatorname{Ind}_{\Gamma_0(p)}^{\Gamma(1)}
 \chi_{\Gamma_0(p)}\right)\\
 &=
 Z_{\Gamma(1)}
 \left(s,\widetilde\rho_\chi\circ\widetilde\pi_p\right).
 \end{aligned}
\end{equation}
For even $\chi$, we have
$\widetilde\rho_\chi=\rho_\chi\circ q_p$ and
$\pi_p=q_p\circ\widetilde\pi_p$, so this becomes
\[
 Z_{\Gamma_0(p)}(s,\chi_{\Gamma_0(p)})
 =
 Z_{\Gamma(1)}(s,\rho_\chi\circ\pi_p).
\]
Taking the logarithmic derivative of the twisted Euler product gives
\[
\frac{Z'_{\Gamma(1)}(s,\widetilde\rho_\chi\circ\widetilde\pi_p)}
     {Z_{\Gamma(1)}(s,\widetilde\rho_\chi\circ\widetilde\pi_p)}
=
\sum_{\gamma\in\mathcal P(\Gamma(1))}
\sum_{j\geq1}
\frac{\ell_\gamma}
     {1-e^{-j\ell_\gamma}}
\Tr\widetilde\rho_\chi(\widetilde\gamma^{\,j})
e^{-js\ell_\gamma},
\qquad
\Re s>1.
\label{eq:twisted-log-derivative-full}
\]
If $\chi$ is even, then
\[
\Tr\widetilde\rho_\chi(\widetilde\gamma^{\,j})
=
\Tr\rho_\chi(\bar\gamma^{\,j}),
\]
and hence the same logarithmic derivative may be written entirely in
the quotient $G_p$.

\subsection{Hashimoto's local factor and the four reduction types}
\label{subsec:hashimoto-local-types}

We first recall the part of Hashimoto’s untwisted formula in \cite{Hashimoto2007} needed below and reinterpret it in terms of reduction modulo \(p\). Let \(\mathcal C(d)\) denote the set of proper equivalence classes of primitive indefinite integral binary quadratic forms of discriminant \(d\). As recalled in Subsection~\ref{subsec:selberg-zeta-class-number}, its \(h(d)\) classes correspond to the primitive hyperbolic conjugacy classes associated with \(d\).
For a representative \(Q=[A,B,C]\), with \(B^2-4AC=d\), let \(\gamma_Q\) denote the corresponding primitive hyperbolic element. For \(j\ge1\), write
\[
\epsilon_d^j=\frac{t_{d,j}+u_{d,j}\sqrt d}{2}.
\]
Then \(\gamma_Q^j\) is represented by
\begin{equation}\label{eq:Q-matrix}
\gamma_{Q,j}=
\begin{pmatrix}
\dfrac{t_{d,j}-Bu_{d,j}}2&-Cu_{d,j}\\[4pt]
Au_{d,j}&\dfrac{t_{d,j}+Bu_{d,j}}2
\end{pmatrix}.
\end{equation}
The Pell relation \(t_{d,j}^2-du_{d,j}^2=4\) gives \(\det\gamma_{Q,j}=1\). Thus its characteristic polynomial is
\begin{equation}\label{eq:char-pol}
X^2-t_{d,j}X+1,
\end{equation}
whose discriminant is \(t_{d,j}^2-4=du_{d,j}^2\).

The results of \cite{Hashimoto2007} identify the multiplicity of the corresponding prime geodesic in the \(\Gamma_0(p)\)-cover with the number of fixed points of
\(\widetilde{\gamma}_{Q,j}\) under the natural left action of
\(\widetilde G_p\) on
\[
\Gamma(1)/\Gamma_0(p)
 \simeq \widetilde G_p/\widetilde B_p
 \simeq \mathbb P^1(\mathbb F_p).
\]
This number is determined by the splitting behavior of the characteristic
polynomial \eqref{eq:char-pol}, together with the
scalar--nonscalar distinction in the repeated-root case.

\begin{proposition}
\label{prop:fixed-points}
Let \(d\in  \Dset\), \(j\ge1\), and let \(Q\) represent a class in \(\mathcal C(d)\).
Let
\(\widetilde{\gamma}_{Q,j}\in\widetilde G_p\) be the reduction modulo \(p\)
of the matrix \(\gamma_{Q,j}\) in \eqref{eq:Q-matrix}. Then the number
of fixed lines of \(\widetilde{\gamma}_{Q,j}\) on
\(\mathbb P^1(\mathbb F_p)\) is independent of \(Q\) and is given by
\begin{equation}
\#\operatorname{Fix}_{\mathbb P^1(\mathbb F_p)}
   \bigl(\widetilde{\gamma}_{Q,j}\bigr)
=
\begin{cases}
p+1,
   & p\mid u_{d,j},\\[2mm]
2,
   & p\nmid u_{d,j},\ 
     \displaystyle\left(\frac{d}{p}\right)=1,\\[2mm]
1,
   & p\nmid u_{d,j},\ p\mid d,\\[2mm]
0,
   & p\nmid u_{d,j},\
     \displaystyle\left(\frac{d}{p}\right)=-1.
\end{cases}
 \label{eq:hashimoto-prime-local-factor}
\end{equation}
Equivalently, the scalar, split, ramified, and inert types have,
respectively, \(p+1\), \(2\), \(1\), and \(0\) fixed lines. 
\end{proposition}

\begin{proof}
If \(p\mid u_{d,j}\), then \eqref{eq:Q-matrix} gives
\[
\widetilde{\gamma}_{Q,j}
   =\frac{t_{d,j}}{2}I=\pm I \qquad \text{in} \quad
\SL_2(\mathbb F_p) ,
\]
so every line in \(\mathbb P^1(\mathbb F_p)\) is fixed.

Suppose now that \(p\nmid u_{d,j}\). If
\(\left(\frac{d}{p}\right)=1\), the characteristic polynomial has two
distinct roots in \(\mathbb F_p\), and hence
\(\widetilde{\gamma}_{Q,j}\) fixes its two eigenlines. If \(p\mid d\),
the characteristic polynomial has a repeated root. The reduction is
non-scalar; otherwise \eqref{eq:Q-matrix}, together with
\(p\nmid u_{d,j}\), would imply \(p\mid A,B,C\), contradicting the
primitivity of \(Q=[A,B,C]\). It therefore has exactly one eigenline.
Finally, if \(\left(\frac{d}{p}\right)=-1\), the characteristic
polynomial is irreducible over \(\mathbb F_p\), so there is no fixed
line.

Thus the fixed-line count depends only on the Pell data
\((t_{d,j},u_{d,j},d)\), and not on \(Q\). Since
\(\operatorname{Ind}_{\widetilde B_p}^{\widetilde G_p}\mathbf 1\) is
the permutation representation on
\[
\widetilde G_p/\widetilde B_p\simeq\mathbb P^1(\mathbb F_p),
\]
its character is precisely this fixed-line count.
\end{proof}

\subsection{Dirichlet-character refinement}
\label{subsec:character-refinement}

We now refine the fixed-line count in Proposition \ref{prop:fixed-points} by recording the
multiplier on each fixed line. 
Let $g\in\widetilde G_p$, and suppose that
\[
   L=x\widetilde B_p
   \in
   \widetilde G_p/\widetilde B_p
   \simeq
   \mathbb P^1(\mathbb F_p)
\]
is fixed by \(g\). Then \(x^{-1}gx\in\widetilde B_p\), so we may
write
\[
   x^{-1}gx
   =
   \begin{pmatrix}
      a & *\\
      0 & a^{-1}
   \end{pmatrix}.
\]
The element \(a\in\mathbb F_p^\times\) depends only on the fixed
line \(L\). 
Replacing $x$ by $xb$ with $b\in\widetilde B_p$ conjugates
$x^{-1}gx$ inside $\widetilde B_p$ and does not change its image in
the abelian quotient $\widetilde B_p/\widetilde U_p$. Hence the
multiplier depends only on $L$.
We call it the multiplier of \(g\) on \(L\) and write
\[
   t(g,L):=a.
\]
By the definition \eqref{eq:induced-rep},
we have
\begin{equation}
   \Tr\widetilde\rho_\chi(g)
   =
   \sum_{\substack{L\in\mathbb P^1(\mathbb F_p)\\gL=L}}
   \chi\bigl(t(g,L)\bigr).
   \label{eq:induced-character-fixed-lines}
\end{equation}
For $d\in \mathcal{D}$ and \(j\geq1\), let
\[
   M_{\Gamma_0(p)}(d,j,\chi)
   :=
   \frac{1}{h(d)}
   \sum_{[\gamma]\in \mathcal C(d)}
   \Tr\widetilde\rho_\chi
   \bigl(\widetilde{\gamma}^{\,j}\bigr).
\] 
Grouping the primitive hyperbolic classes in the logarithmic derivative
of the Selberg zeta function according to their discriminants gives
\begin{equation}\label{eq:zeta-log-derivative-class-number}
\begin{aligned}
\frac{Z_{\Gamma_0(p)}'(s,\chi_{\Gamma_0(p)})}
     {Z_{\Gamma_0(p)}(s,\chi_{\Gamma_0(p)})}
=&
\frac{Z_{\Gamma(1)}'(s,\widetilde{\rho}_{\chi}\circ \widetilde{\pi}_p)}
     {Z_{\Gamma(1)}(s,\widetilde{\rho}_{\chi}\circ \widetilde{\pi}_p)}\\
=&
\sum_{d\in \Dset} h(d)
\sum_{j\geq 1}
\frac{2\log\epsilon_d}
     {1-\epsilon_d^{-2j}}\,
M_{\Gamma_0(p)}(d,j,\chi)\,
\epsilon_d^{-2js},
\qquad \Re s>1.
\end{aligned}
\end{equation}
The four reduction types determine this coefficient $M_{\Gamma_0(p)}(d,j,\chi)$ explicitly.

\begin{proposition}
\label{prop:character-local-formula}
For an odd prime $p$, $d\in \Dset$, a Dirichlet character $\chi$ modulo $p$, and $j\geq1$,
\begin{equation}
   M_{\Gamma_0(p)}(d,j,\chi)
   =
   \begin{cases}
      (p+1)\chi(t_{d,j}/2),
         & p\mid u_{d,j},\\[1mm]
      \chi(\alpha_+)+\chi(\alpha_-),
         & p\nmid u_{d,j},\
           \displaystyle\left(\frac{d}{p}\right)=1,\\[2mm]
      \chi(t_{d,j}/2),
         & p\nmid u_{d,j},\ p\mid d,\\[1mm]
      0,
         & p\nmid u_{d,j},\
           \displaystyle\left(\frac{d}{p}\right)=-1,
   \end{cases}
   \label{eq:character-local-four-cases}
\end{equation}
where, in the split case, if \(r^2\equiv d\pmod p\), then
\[
   \alpha_\pm
   :=
   \frac{t_{d,j}\pm u_{d,j}r}{2}
   \in\mathbb F_p^\times,
   \qquad
   \alpha_+\alpha_-=1.
\]
\end{proposition}

\begin{proof}
In the scalar case,
\[
   \widetilde{\gamma}^{\,j}
   =
   \frac{t_{d,j}}{2}I \qquad \text{in} \quad
\SL_2(\mathbb F_p) 
\]
fixes all \(p+1\) lines, with the same multiplier
\(t_{d,j}/2\) on each line. Hence
\[
   \Tr\widetilde\rho_\chi
   \bigl(\widetilde{\gamma}^{\,j}\bigr)
   =
   (p+1)\chi(t_{d,j}/2).
\]
In the split case, \(\widetilde{\gamma}^{\,j}\) has two fixed
eigenlines with eigenvalues \(\alpha_+\) and \(\alpha_-\).
Equation \eqref{eq:induced-character-fixed-lines} therefore gives
\[
   \Tr\widetilde\rho_\chi
   \bigl(\widetilde{\gamma}^{\,j}\bigr)
   =
   \chi(\alpha_+)+\chi(\alpha_-).
\]
In the ramified case, the reduction is non-scalar with a unique
fixed line and repeated eigenvalue \(t_{d,j}/2\), giving
\[
   \Tr\widetilde\rho_\chi
   \bigl(\widetilde{\gamma}^{\,j}\bigr)
   =
   \chi(t_{d,j}/2).
\]
Finally, an inert element has no fixed line over
\(\mathbb F_p\), so its trace is zero. 
These values depend only
on the Pell data and are therefore the same for every class in
\(\mathcal{C}(d)\). Averaging over \(|\mathcal{C}(d)|=h(d)\) proves the result.
\end{proof}

The result in
\eqref{eq:hashimoto-prime-local-factor} is recovered from
\eqref{eq:character-local-four-cases} by taking
\(\chi=\mathbf 1\). The character-weighted formula is therefore a
refinement of the untwisted result.

If \(\chi\) is even, then
\(\widetilde\rho_\chi\) is trivial on \(\{\pm I\}\) and descends
to
\[
   \rho_\chi
   :=
   \operatorname{Ind}_{B_p}^{G_p}\chi,
   \qquad
   G_p=\PSL_2(\mathbb F_p),
   \qquad
   H_p=B_p/U_p\simeq\mathbb F_p^\times/\{\pm1\}.
\]
Writing
\[
   a_\pm=[\alpha_\pm]\in H_p,
   \qquad
   a_-=a_+^{-1},
\]
the even-character formula becomes
\[
   M_{\Gamma_0(p)}(d,j,\chi)
   =
   \begin{cases}
      (p+1)\chi(1),
         & p\mid u_{d,j},\\[1mm]
      \chi(a_+)+\chi(a_-),
         & p\nmid u_{d,j},\
           \displaystyle\left(\frac{d}{p}\right)=1,\\[2mm]
      \chi(1),
         & p\nmid u_{d,j},\ p\mid d,\\[1mm]
      0,
         & p\nmid u_{d,j},\
           \displaystyle\left(\frac{d}{p}\right)=-1.
   \end{cases}
   \label{eq:even-character-local-four-cases}
\]

%\beb Thus the four local types, together with the projective multiplier in the split case and the central sign for odd characters, determine all the character-weighted local coefficients. \eb \ber [Is the meaning of this sentence clear? Do we need this?] [This is not needed here. Actually this summarizes the content of the following corollary. We can remove this.]\er

\begin{corollary}\label{cor:character-summed-weights}
For $j=1$, summing the coefficient in Proposition~\ref{prop:character-local-formula}
over the nontrivial even characters gives
\begin{equation}
\label{eq:even-local-type-weight}
\sum_{\substack{\chi\bmod p\\ \chi(-1)=1,\ \chi\neq \mathbf 1}}
M_{\Gamma_0(p)}(d,1,\chi)
=
\begin{cases}
\dfrac{(p+1)(p-3)}{2}, & p\mid u_d,\\[4pt]
-2, & p\nmid u_d,\ \left(\dfrac d p\right)=1,\\[4pt]
\dfrac{p-3}{2}, & p\nmid u_d,\ p\mid d,\\[4pt]
0, & p\nmid u_d,\ \left(\dfrac d p\right)=-1.
\end{cases}
\end{equation}
For the odd characters, one has
\begin{equation}
\label{eq:odd-local-type-weight} 
\sum_{\substack{\chi\bmod p\\ \chi(-1)=-1}}
M_{\Gamma_0(p)}(d,1,\chi)
=
\begin{cases}
\sigma_p(d)\dfrac{(p+1)(p-1)}{2}, & p\mid u_d,\\[4pt]
0, & p\nmid u_d,\ \left(\dfrac d p\right)=1,\\[4pt]
\sigma_p(d)\dfrac{p-1}{2}, & p\nmid u_d,\ p\mid d,\\[4pt]
0, & p\nmid u_d,\ \left(\dfrac d p\right)=-1,
\end{cases}
\end{equation}
where 
\begin{equation*}
\sigma_p(d)\in\{\pm1\},
\qquad
\frac{t_d}{2}\equiv \sigma_p(d)\pmod p.
\label{eq:sigma-p-definition}
\end{equation*}
\end{corollary}
\begin{proof}
For even characters, character orthogonality on
\[
H_p=\mathbb F_p^\times/\{\pm1\}
\]
gives
\begin{equation*}
\sum_{\substack{\chi\in\widehat H_p\\\chi\neq \mathbf 1}}
\chi(a)
=
\begin{cases}
\dfrac{p-3}{2},
&
a=1,
\\[2mm]
-1,
&
a\neq1.
\end{cases}
\label{eq:even-character-orthogonality}
\end{equation*}
In the split case, 
where the reduction has two distinct eigenvalues,
neither multiplier class $a_+$ nor $a_-$
is the identity of $H_p$. Indeed, $a_\pm=1$ would imply
$\alpha_+ \in \{\pm1\}$. Since $\alpha_-=\alpha_+^{-1}$, this would imply
$\alpha_+=\alpha_-$, contradicting the distinctness of the two eigenvalues. Therefore, each of the two fixed lines contributes $-1$
to the sum over nontrivial even characters, giving $-2$ for the split case in \eqref{eq:even-local-type-weight}.
For the odd characters, ordinary character orthogonality on
$\mathbb F_p^\times$ gives
\[
\sum_{\substack{\chi\bmod p\\\chi(-1)=-1}}\chi(y)
=
\begin{cases}
\dfrac{p-1}{2},
&
y=1,
\\[2mm]
-\dfrac{p-1}{2},
&
y=-1,
\\[2mm]
0,
&
y\notin\{\pm1\}.
\end{cases}
\]
Combining these and \eqref{eq:character-local-four-cases} gives the claimed equalities.
\end{proof}

For each $d\in \Dset$, define the character-summed weights
\begin{equation}\label{eq:character-summed-weights}
\begin{aligned}
w_p^+(d)
&:=
\sum_{\substack{\chi\bmod p\\
\chi(-1)=1,\ \chi\neq\mathbf 1}}
M_{\Gamma_0(p)}(d,1,\chi),\\
w_p^-(d)
&:=
\sum_{\substack{\chi\bmod p\\
\chi(-1)=-1}}
M_{\Gamma_0(p)}(d,1,\chi).
\end{aligned}
\end{equation}
Corollary~\ref{cor:character-summed-weights} shows that the resulting sums over nontrivial even or odd characters can be evaluated without enumerating the Dirichlet characters or individual quadratic-form classes.
For each discriminant $d$, it is enough to compute $(t_d,u_d,d)$ modulo $p$, determine its reduction type, and, for the odd sum, determine the sign $\sigma_p(d)$.
The weights are explicitly given by
\begin{equation}\label{eq:weights-by-reduction-type}
\begin{array}{c|cc}
\text{Reduction type of }d
& w_p^+(d) & w_p^-(d)\\
\hline
\text{scalar}
& \dfrac{(p+1)(p-3)}{2}
& \dfrac{(p+1)(p-1)}{2}\sigma_p(d)\\[6pt]
\text{split}
& -2 & 0\\[6pt]
\text{ramified}
& \dfrac{p-3}{2}
& \dfrac{p-1}{2}\sigma_p(d)
\end{array}
\end{equation}
These formulas form the basis of the algorithms used in Sections~\ref{sec:experimental-setup} and~\ref{sec:local-type-contributions}.

\section{Character-twisted sums of class numbers and numerical experiments}
\label{sec:experimental-setup}

In this section, we use the formulas of Section~\ref{sec:class-number-geodesics} to define the
character-weighted class number sums that will be studied in the subsequent sections. After
determining the natural density normalization, we present 
numerical data and examine the resulting oscillatory behavior for fixed primes. These
experiments provide the initial evidence for Conjecture~\ref{conj:ap-experimental}.

\subsection{Arithmetic sums and density normalization}

We first specify the discriminant families naturally associated with the
congruence subgroups \(\Gamma(p)\) and \(\Gamma_0(p)\). For
\(d\in \Dset\) with the Pell data $(t_d, u_d)$
satisfying \eqref{eq:Pell-data},
define
\begin{equation*}
 \Dset_p
:=
\{d\in \Dset:p\mid u_d\}.
\label{eq:Dp}
\end{equation*}
For \(d\in \Dset_p\), the corresponding primitive hyperbolic
elements reduce to \(\pm I\) modulo \(p\) and are therefore of scalar
type. 
Their reductions modulo \(p\) are \(\pm I\); equivalently, up to multiplication by \(-I\), they lie in the principal congruence subgroup \(\Gamma(p)\).
Thus \( \Dset_p\) is precisely the scalar-type
discriminant family naturally associated with $\Gamma(p)\cup -\Gamma(p)$.
The family associated with \(\Gamma_0(p)\) is
\begin{equation*}
 \Dset_{0,p}
:=
\left\{
d\in \Dset: \quad
p\mid u_d
\quad \text{or}\quad
\left(\frac{d}{p}\right)\in\{0,1\}
\right\}.
\label{eq:D0p}
\end{equation*}
Indeed, the Pell equation in \eqref{eq:Pell-data}, together with Proposition~\ref{prop:fixed-points}, shows that
\(d\in \Dset_{0,p}\) precisely when the reduced hyperbolic element
has a fixed line in \(\mathbb P^1(\mathbb F_p)\). Since 
\(\widetilde B_p\) is the stabilizer of a line in \(\mathbb P^1(\mathbb F_p)\), the reduced element has a fixed line if and only if it is conjugate in \(\widetilde G_p\) to an element of \(\widetilde B_p\). 
This is precisely the condition relevant to the passage from \(\Gamma(1)\) to \(\Gamma_0(p)\). Consequently, \( \Dset_{0,p}\) consists exactly of the scalar, split, and ramified reduction types; the inert type, which has no fixed line over \(\mathbb F_p\), does not occur.

For real \(x\geq1\), let
\[
 \Dset(x):=\{d\in \Dset:\epsilon_d\leq x\}.
\]
We then define the corresponding truncated families by
\begin{equation}
 \Dset_p(x):= \Dset_p\cap \Dset(x),
\qquad
 \Dset_{0,p}(x):= \Dset_{0,p}\cap \Dset(x).
\label{eq:truncated-discriminant-families}
\end{equation}
Thus
\[
 \Dset_{0,p}(x)
=
 \Dset_p(x)
\sqcup
\left\{
d\in \Dset(x):
p\nmid u_d,\ 
\left(\frac{d}{p}\right)\in\{0,1\}
\right\}.
\]
For each of the three noninert reduction types (scalar, split, and ramified),
we define
\begin{equation*}\label{eq:D0p-type-families}
 \Dset_{0,p}^{T}(x)
:=
\{d\in \Dset_{0,p}(x):
d\text{ has reduction type }T\},
\end{equation*}
where the reduction types are determined by the Pell conditions in \eqref{eq:types}. Thus
\[
 \Dset_{0,p}(x)
=
 \Dset_{0,p}^{\Scalar}(x)
\sqcup
 \Dset_{0,p}^{\Split}(x)
\sqcup
 \Dset_{0,p}^{\Ramified}(x),
\]
with
\[
 \Dset_{0,p}^{\Scalar}(x)
= \Dset_p(x).
\]

To determine the natural density factors used in the normalization, we need the following:

\begin{lemma}\label{lem:count}
For every odd prime $p$,
\[
\sum_{a\in\mathbb F_p}
\left(\frac{a^2-4}{p}\right)=-1.
\]
\end{lemma}

\begin{proof} For each \(a\in\mathbb F_p\), the number of \(y\in\mathbb F_p\) satisfying $y^2=a^2-4$ is
$1+\left(\frac{a^2-4}{p}\right)$.
Indeed, this number is \(2\), \(1\), or \(0\), according as \(a^2-4\) is a nonzero square, zero, or a nonsquare. Therefore, the total number \(N\) of solutions \((a,y)\in\mathbb F_p^2\) is
$$ N=\sum_{a\in\mathbb F_p} \left(1+\left(\frac{a^2-4}{p}\right)\right) =p+\sum_{a\in\mathbb F_p}\left(\frac{a^2-4}{p}\right). $$
It remains to compute \(N\). The equation \(y^2=a^2-4\) can be rewritten as
$(a-y)(a+y)=4$.
Since the right-hand side is nonzero, \(r:=a-y\) belongs to \(\mathbb F_p^\times\), and then
$a+y=4r^{-1}$.
Conversely, because \(p\) is odd, every \(r\in\mathbb F_p^\times\) determines a unique solution

$$ a=\frac{r+4r^{-1}}2, \qquad y=\frac{4r^{-1}-r}2. $$
Thus the solution set is in bijection with \(\mathbb F_p^\times\), and hence \(N=p-1\). Comparing the two expressions for \(N\) gives the claimed equality.
\end{proof}

\begin{proposition}\label{prop:density}
For every odd prime \(p\), 
assume that, 
within the subfamily \(p\nmid u_d\), each residue class of \(t_d\) modulo \(p\) occurs with asymptotic proportion \(1/p\), counting each discriminant once.
Then
\begin{equation*}\label{eq:normalizing-constant}
  |\Dset_{0,p}(x)|\sim C_p x,
  \qquad
  C_p=
  \frac{35}{32}
\,  \frac{p^3+p^2+3p-1}{p(p^2+1)}.
\end{equation*}
Equivalently,
\[
  \frac{|\Dset_{0,p}(x)|}{|\Dset(x)|}
  \longrightarrow
  \frac{p^3+p^2+3p-1}{2p(p^2+1)}.
\]
\end{proposition}

\begin{proof}
The original class number counting result in \cite[Proposition~4.1]{Sarnak1982} gives
\[
  |\Dset(x)|\sim \frac{35}{16}x,
\]
and the scalar subfamily satisfies
\[
  |\Dset_p(x)|\sim \frac{35}{8(p^2+1)}x.
\]
Hence the scalar proportion is \(2/(p^2+1)\), and the complementary set
\[ \Dset_p^{\mathrm{comp}}(x)
:= \Dset(x)\setminus \Dset_p(x)=\{d\in \Dset(x):p\nmid u_d\}
\]
has proportion \((p^2-1)/(p^2+1)\).

For $d\in \Dset_p^{\mathrm{comp}}(x)$, the Pell equation gives
\[
  d\equiv (t_d^2-4)u_d^{-2}\pmod p.
\]
Since $p\nmid u_d$, the element $u_d^{-1}\in\mathbb F_p^\times$
is defined and $u_d^{-2}=(u_d^{-1})^2$ is a nonzero square.
It follows that
\[
\left(\frac{d}{p}\right)
=
\left(\frac{t_d^2-4}{p}\right)
\left(\frac{u_d^{-2}}{p}\right)
=
\left(\frac{t_d^2-4}{p}\right).
\]
Under the stated local equidistribution hypothesis,
$\Leg d p$ has the same limiting distribution as
$\left(\frac{a^2-4}{p}\right)$ with $a\in\F_p$.
Moreover, $\left(\frac{a^2-4}{p}\right)=0$ precisely when $a=\pm2$.
If $N_+$ and $N_-$ denote the numbers of $a\in\mathbb F_p$
for which this symbol equals $1$ and $-1$, respectively,
then Lemma~\ref{lem:count} gives
\[
N_++N_-=p-2,
\qquad
N_+-N_-=\sum_{a\in\mathbb F_p}
\left(\frac{a^2-4}{p}\right)=-1.
\]
These identities yield
\[
N_+=\frac{p-3}{2}
\qquad \text{and} \qquad
N_-=\frac{p-1}{2}.
\]
Thus the values $1,0,-1$ occur, respectively, for $(p-3)/2$, $2$, and $(p-1)/2$ residue classes $a\in\mathbb F_p$, corresponding to the split, ramified, and inert reduction types.
% Consequently, the conditional noninert proportion in \( \Dset_p^{\mathrm{comp}}(x)\) is
% \[
%   \frac{p-3}{2p}+\frac{2}{p}=\frac{p+1}{2p}.
% \]
% Since $\Dset_{0,p}(x)$ is the disjoint union of $\Dset_p(x)$ and the noninert discriminants in $\Dset_p^{\mathrm{comp}}(x)$, we obtain
% \[
%   \frac{|\Dset_{0,p}(x)|}{|\Dset(x)|}
%   \longrightarrow
%   \frac{2}{p^2+1}
%   +
%   \frac{p^2-1}{p^2+1}\cdot\frac{p+1}{2p}
%   =
%   \frac{p^3+p^2+3p-1}{2p(p^2+1)}.
% \]
Combining these counts with the proportions of $\Dset_p(x)$ and $\Dset_p^{\mathrm{comp}}(x)$ in $\Dset(x)$ obtained above, we find that the scalar, split, and ramified types have limiting proportions given, respectively, by
\begin{equation}\label{eq:unweighted-type-densities}
\frac{2}{p^2+1},\qquad
\frac{p^2-1}{p^2+1}\cdot \frac{p-3}{2p},\qquad
\frac{p^2-1}{p^2+1}\cdot\frac{2}{p}.
\end{equation}
Since $\Dset_{0,p}(x)$ consists exactly of these three types,
we obtain
\[
  \frac{|\Dset_{0,p}(x)|}{|\Dset(x)|}
  \longrightarrow
  \frac{2}{p^2+1}
  +\frac{(p^2-1)(p-3)}{2p(p^2+1)}
  +\frac{2(p^2-1)}{p(p^2+1)}
  =
  \frac{p^3+p^2+3p-1}{2p(p^2+1)}.
\]
Multiplying by $|\Dset(x)|\sim(35/16)x$ gives
the stated value of $C_p$.

\end{proof}

\begin{remark}
The constant $C_p$ is used only as a density normalization for comparing different primes.
None of the other mathematical results in this paper depends on this equidistribution assumption.
\end{remark}

Using the character-summed weights $w_p^\pm(d)$ defined in \eqref{eq:character-summed-weights}, we set, for $d\in \Dset$,
\[
a_p^\pm(d):=h(d)\,w_p^\pm(d).
\]
Recalling \eqref{eq:truncated-discriminant-families}, for real $x\geq1$, define the even and odd character-twisted sums of class numbers by
\begin{equation*}
A_p^\pm(x)
:=
\sum_{d\in \Dset_{0,p}(x)} a_p^\pm(d).
\label{eq:Ap-pm-local-summands}
\end{equation*}
For each of the three noninert reduction types $T$ (scalar, split, or ramified), we also define
\begin{equation}\label{eq:Ap-type-components}
A_p^{T,\pm}(x)
:=
\sum_{d\in \Dset_{0,p}^{T}(x)}
a_p^\pm(d).
\end{equation}
Since the three reduction types partition
$ \Dset_{0,p}(x)$, we have
\begin{equation}\label{eq:Ap-type-decomposition}
A_p^\pm(x)
=
A_p^{\Scalar,\pm}(x)
+
A_p^{\Split,\pm}(x)
+
A_p^{\Ramified,\pm}(x).
\end{equation}
We further put
\[
A_p^T(x):=A_p^{T,+}(x)+A_p^{T,-}(x),
\qquad
A_p(x):=A_p^+(x)+A_p^-(x).
\]

For a Dirichlet character $\chi$ modulo $p$, we also define
\begin{equation}
A_{p,\chi}(x)
:=
\sum_{d\in \Dset_{0,p}(x)}
h(d)\,
M_{\Gamma_0(p)}(d,1,\chi)
\label{eq:def-A-p-chi}
\end{equation}
so that
\begin{equation*}
A_p^+(x)
=
\sum_{\substack{\chi\bmod p\\
                 \chi(-1)=1,\ \chi\neq \mathbf 1}}
A_{p,\chi}(x),
\qquad
A_p^-(x)
=
\sum_{\substack{\chi\bmod p\\
                 \chi(-1)=-1}}
A_{p,\chi}(x).
\label{eq:Ap-plus-minus-definition}
\end{equation*}

\subsection{\texorpdfstring
  {Experiments for \(A_p^+(x)\), \(A_p^-(x)\), and \(A_p(x)\)}
  {Experiments for A-plus(x), A-minus(x), and A(x)}}

The scalar, split, and ramified discriminants contributing to
\(A_p^\pm(x)\) occur in \( \Dset(x)\) with different asymptotic
frequencies and carry different character-summed weights.

By \eqref{eq:weights-by-reduction-type}, the even-character weight \(w_p^+(d)\) takes the values
\[
\frac{(p+1)(p-3)}{2},\qquad
-2,\qquad
\frac{p-3}{2},
\]
for the scalar, split, and ramified types, respectively.
Under the assumption of Proposition~\ref{prop:density}, their densities in $ \Dset$ are given by \eqref{eq:unweighted-type-densities}.
Thus the scalar, split, and ramified types have densities of orders \(p^{-2}\), \(1\), and \(p^{-1}\), but weights of orders \(p^2\), \(1\), and \(p\). 
The contributions summed over the respective reduction types may therefore be comparable,
while the negative split weight permits substantial cancellation with
the positive scalar and ramified weights.

A similar balance occurs in $A_p^-(x)$. 
By \eqref{eq:weights-by-reduction-type}, the odd-character weight \(w_p^-(d)\) takes the values
\[
\sigma_p(d)\frac{(p+1)(p-1)}{2},\qquad
0,\qquad
\sigma_p(d)\frac{p-1}{2},
\]
for the scalar, split, and ramified types, respectively.
Thus the scalar and ramified weights have orders $p^2$ and $p$, respectively, and carry the additional sign $\sigma_p(d)$, while the split contribution vanishes.

Consequently, the combined sum
\[
A_p(x)=A_p^+(x)+A_p^-(x)
\]
may exhibit further cancellation between the even and odd parts.

\begin{figure}[p]
\centering
\includegraphics[
  width=0.95\textwidth,
  height=0.90\textheight,
  keepaspectratio
]{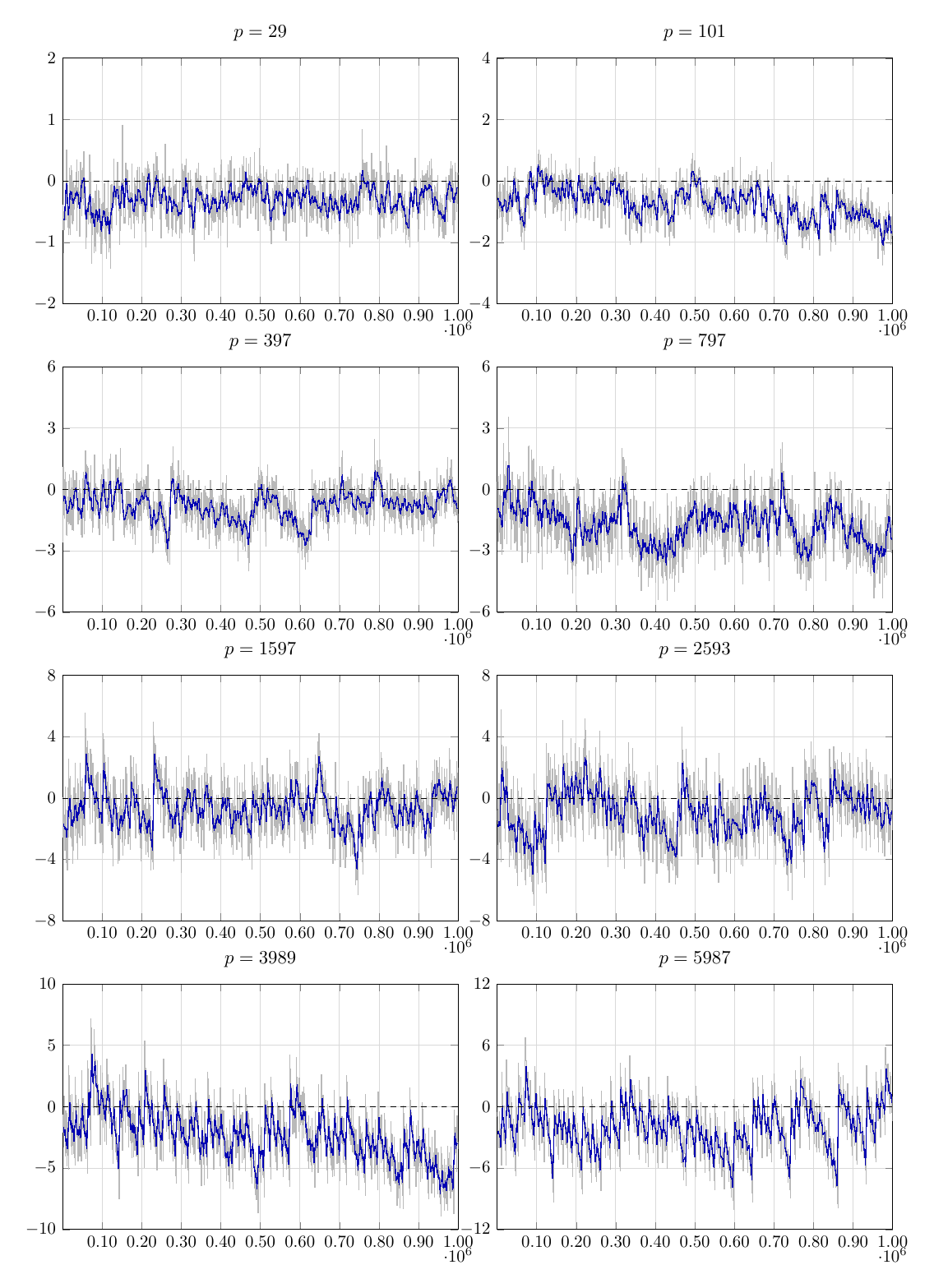}
\caption{\sf \small The ratios $A_p^+(x)/(\sqrt p\,C_p x)$ for selected primes $p$. \ApmMACaption}
\label{fig:ap-plus-over-cpx}
\end{figure}

\begin{figure}[p]
\centering
\includegraphics[
  width=0.95\textwidth,
  height=0.90\textheight,
  keepaspectratio
]{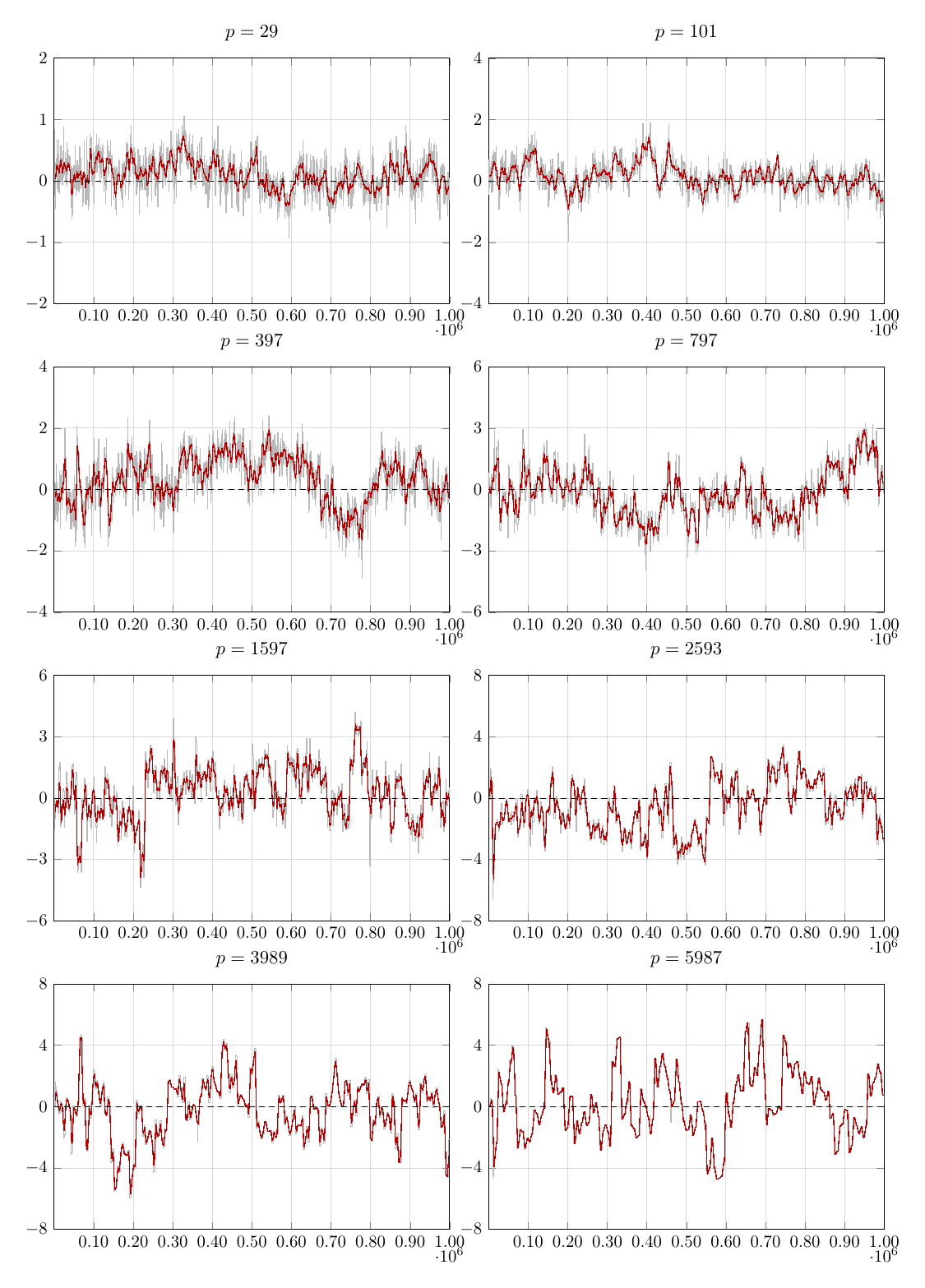}
\caption{\sf \small The ratios $A_p^-(x)/(\sqrt p\,C_p x)$ for selected primes $p$. \ApmMACaption}
\label{fig:ap-minus-over-cpx}
\end{figure}

\begin{figure}[p]
\centering
\includegraphics[
  width=0.95\textwidth,
  height=0.90\textheight,
  keepaspectratio
]{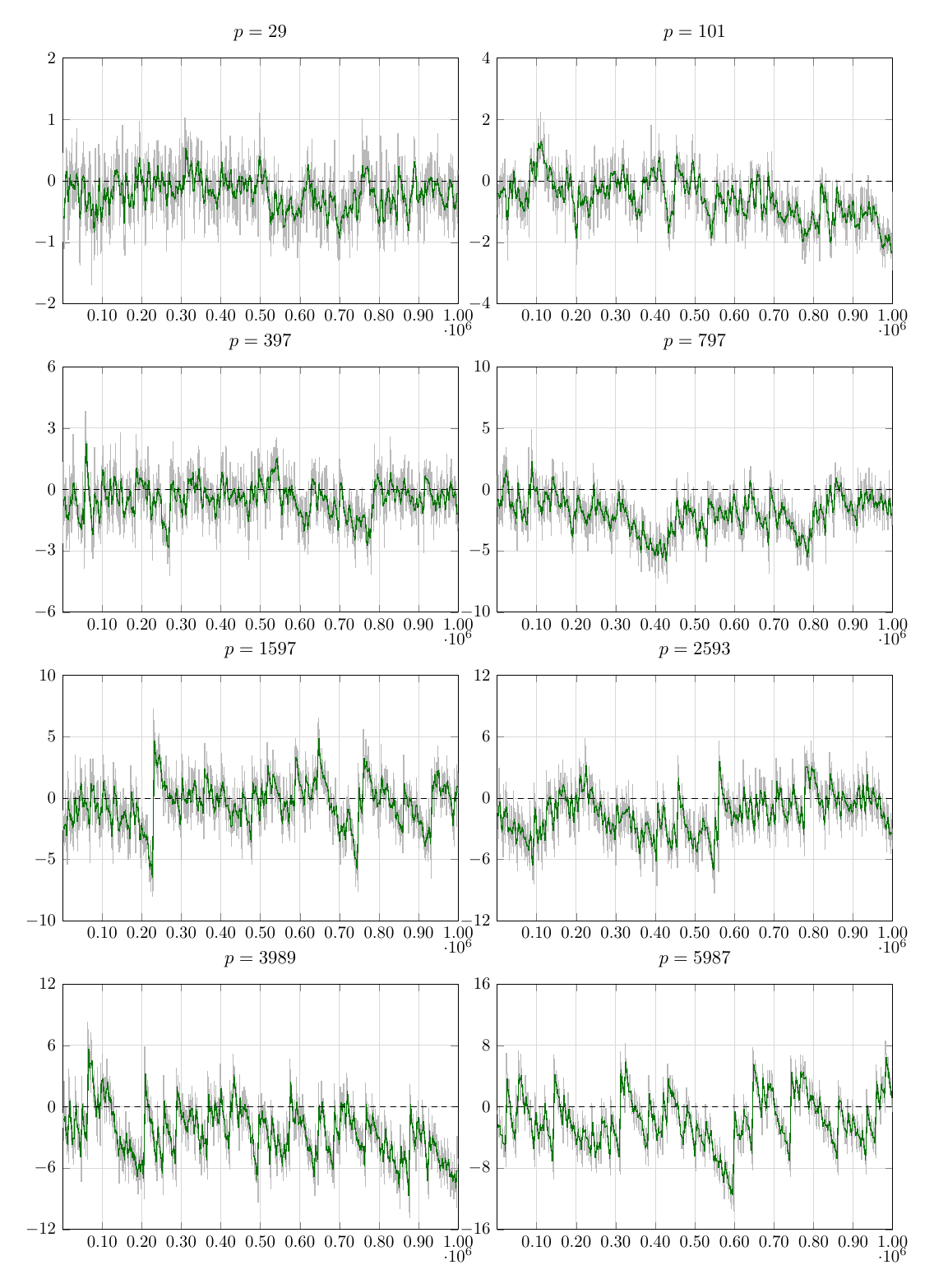}
\caption{\sf \small The ratios $A_p(x)/(\sqrt p\,C_p x)$ for selected primes $p$. \ApmMACaption}
\label{fig:ap-sum-over-cpx}
\end{figure}

Figures~\ref{fig:ap-plus-over-cpx}--\ref{fig:ap-sum-over-cpx}
show three experimental panels for
\begin{equation}
\label{eq:primes}    
  p=29,101,397,797,1597,2593,3989,5987.
\end{equation}
They plot
\[
\frac{A_p^+(x)}{\sqrt p\, C_p x },\qquad
\frac{A_p^-(x)}{\sqrt p \, C_p x},\qquad
\frac{A_p(x)}{\sqrt p\,C_p x}=
  \frac{A_p^+(x)+A_p^-(x)}{\sqrt p\, C_p x}.
\]
The horizontal axis is \(x\), sampled up to \(10^6\).  Raw values are plotted in gray every
\(500\) in \(x\), and the colored curve averages the normalized values at all integer points in
\([x-2500,x+2500]\).  These panels are the basic evidence for the normalization by \(x\sqrt p\).

The important feature of Figures~\ref{fig:ap-plus-over-cpx}--\ref{fig:ap-sum-over-cpx}
is not merely the oscillation of the curves, but their scale.  For the primes in \eqref{eq:primes}, both even
and odd sums remain of size comparable with \(x\sqrt p\) through the plotted range. However, the observation made in Remark \ref{re:Arithmetic jumps and the conjectural bound} corrects this bound and suggests the estimate
\begin{equation*}
\label{eq:odd-conjecture-A^-_p}
A_p^{\pm}(x)\ll_{\epsilon}
 (px)^{1+\epsilon}.
\end{equation*}
Hence the combined statistic
\(A_p(x)=A_p^+(x)+A_p^-(x)\) is expected to satisfy a bound of the same
order, even though its two components arise from sums over characters of different parity.
As we explain in Section \ref{sec:dirichlet-spectral-gap}, in the nontrivial even Dirichlet sector, such estimates retain the
\(px\)-scale after insertion of the primitive orbital weight, with
at most one additional logarithmic factor.  Passage
to the full trace formula also requires control of other geometric and scattering terms, as discussed in
Section~\ref{sec:dirichlet-spectral-gap}.

\FloatBarrier

\begin{remark} \label{re:Arithmetic jumps and the conjectural bound} 
From the plots in Figures \ref{fig:ap-plus-over-cpx} and \ref{fig:ap-minus-over-cpx}, we may expect 
\begin{equation}\label{eq:expected-inequality}
|A_p^\pm(x)|\ll_\epsilon p^{1/2}x\, (px)^{\epsilon}.
\end{equation}
However,  
the gray values are sampled every \(500\) units in \(x\), while the colored curves average the normalized values over all integer points in \([x-2500,x+2500]\), an interval of length \(5000\). These samplings can conceal the underlying jumps and make the colored curves smoother by averaging over larger intervals.

To explain this, for each trace \(t\), define
\[
\Delta A_p^\pm(\alpha_t)
:=A_p^\pm(\alpha_t)-A_p^\pm(\alpha_t^-)
\]
where \[
A_p^\pm(\alpha_t^-)
:=\lim_{x\uparrow\alpha_t}A_p^\pm(x).
\]
We also recall the following two classical results. First,
we have the class number formula
\[
h(d)\log\epsilon_d=\sqrt{d}\,L(1,\chi_d),
\qquad
\chi_d(n)=\left(\frac{d}{n}\right)
\] 
where $L(s,\chi_d)$ is the Dirichlet L-function.
The inequality from Siegel’s theorem in \cite{Siegel1935} states 
\[
L(1,\chi_d)\gg_\epsilon d^{-\epsilon}
\qquad\text{for every }\varepsilon>0.
\]
Hence, combining these we have
\begin{equation}\label{eq:class number-lower bound}
h(d)\gg_\epsilon\frac{d^{1/2-\epsilon}}{\log\epsilon_d}.
\end{equation}

Now, let \(\alpha_t=(t+\sqrt{t^2-4})/2\) be the larger eigenvalue of a hyperbolic matrix with positive trace \(t>2\), so \(\alpha_{t_d}=\epsilon_d\) is precisely the cutoff at which the classes of discriminant \(d\) enter the sum.
At a trace \(t\equiv2\pmod p\), the ramified weights in \(A_p^\pm\) are \(c_p^+=(p-3)/2\) and \(c_p^-=(p-1)/2\), while the scalar weights are \((p+1)c_p^\pm\) by \eqref{eq:weights-by-reduction-type}. Thus a class number \(h(d)\) produces a positive jump multiplied by a weight of order \(p\) or \(p^2\).

The discriminants \(p(p+4)\) and \(p^2+4\) have ramified and scalar reduction types, and have fundamental norm one solutions \((t,u)=(p+2,1)\) and \((p^2+2,p)\), respectively. Both discriminants have size \(\asymp p^2\), and \(\log\epsilon_d\asymp\log p\).
Multiplying the lower bounds in \eqref{eq:class number-lower bound} for these two class numbers by the corresponding ramified and scalar weights yields the following estimates for the jumps:
\[
\Delta A_p^\pm(\alpha_{p+2})
 \gg_\epsilon\frac{p^{2-2\epsilon}}{\log p},
\qquad
\Delta A_p^\pm(\alpha_{p^2+2})
 \gg_\epsilon\frac{p^{3-2\epsilon}}{\log p}.
\]
A jump \(J\) forces one adjacent value to have magnitude at least \(J/2\). Since the corresponding cutoffs satisfy \(x\asymp p\) and \(x\asymp p^2\) respectively, these jumps are compatible with
\[
|A_p^\pm(x)|\ll_\epsilon (px)^{1+\epsilon}
\]
rather than \eqref{eq:expected-inequality}. 
Hence, nearby opposite-sign jumps can largely compensate each other, so their effect on the sum persists only over a short interval. Averaging over an interval of length \(5000\) reduces this effect, allowing the colored curve to remain smooth despite large local jumps.
In the even sum, negative split contributions before or after a positive jump may compensate it between sampled cutoffs. In the odd sum, whose split weight is zero, compensation instead involves opposite-sign ramified or scalar contributions.

\end{remark}

\begin{remark}\label{re:Paired cancellation in the odd sum}

For the larger primes, the gray values in Figure~\ref{fig:ap-minus-over-cpx} follow the red moving average curve more closely than the corresponding gray values follow the colored curve in Figure~\ref{fig:ap-plus-over-cpx}. To understand this difference, we examine the signs and spacing of the individual jumps.

Put
\[
R_p(t)=\sum_{\substack{d:\,t_d=t\\p\nmid u_d,\;p\mid d}}h(d),
\qquad
S_p(t)=\sum_{\substack{d:\,t_d=t\\p\mid u_d}}h(d),
\qquad
b_p(t)=R_p(t)+(p+1)S_p(t).
\]
For $p\ge5$ and $k\ge1$, put $t_\pm=kp\pm2$ and
$b_\pm=b_p(t_\pm)$.
At \(t_\pm=kp\pm2\), only ramified and scalar discriminants can contribute, and \(b_p(t_\pm)=R_p(t_\pm)+(p+1)S_p(t_\pm)\) combines them with their relative character weights to determine the entire jump. Either mass may be zero.
Then, \eqref{eq:weights-by-reduction-type} gives
\[
\Delta A_p^+(\alpha_{t_\pm})=\frac{p-3}{2}b_\pm,
\qquad
\Delta A_p^-(\alpha_{t_\pm})=\pm\frac{p-1}{2}b_\pm.
\]
The even jumps therefore reinforce each other, whereas the odd jumps have opposite signs.
When \(b_+\) and \(b_-\) are close, the second jump almost compensates the first, returning the sum close to its preceding value. These jumps happen within the interval of size
\[
x_+-x_-:=\alpha_{t_+}-\alpha_{t_-}=4+O((kp)^{-2}),
\]
while the centers of successive pairs are separated by
approximately $p$. 
Such jumps within narrow intervals become increasingly sparse as \(p\) grows. Sampling at intervals of \(500\) may conceal them, while averaging over \([x-2500,x+2500]\) substantially reduces their amplitude in the colored curve. Thus the smoother appearance of the plotted curves can mask large local fluctuations. In the odd case, nearby ramified contributions have opposite signs and can cancel; in the even case, the positive jumps must instead be offset by negative split contributions at other traces.
This explains why the colored curves in the odd case remain closer to zero than those in the even case.
\end{remark}

\section{Asymptotics by reduction type}
\label{sec:cancellation-thm}

The experiments of Section~\ref{sec:experimental-setup} suggest substantial cancellation among
the scalar, split, and ramified contributions to \(A_p^\pm(x)\).
In this section, we establish asymptotic formulas for the individual scalar, split, and ramified contributions to $A_p^+(x)$. As a consequence, we show that their leading terms cancel exactly.

\medskip

We first provide a {Chebotarev-type geodesic formula}.

\begin{lemma}
\label{lem:finite-quotient-pgt}
Fix an odd prime $p$, and let $f$ be a class function on
$G_p=\mathrm{PSL}_2(\mathbb F_p)$.
Then, as $X\to\infty$,
\[
\sum_{\substack{\gamma\in \mathcal P(\Gamma(1))\\N(\gamma)\leq X}}
 f\bigl(\pi_p(\gamma)\bigr)
 =
 \langle f,\mathbf 1\rangle_{G_p}\operatorname{Li}(X)
 +O_{p,f,\epsilon}( X^{25/36+\epsilon})    
\]
where
\[
\langle f_1,f_2\rangle_{G_p}
 :=
 \frac{1}{|G_p|}
 \sum_{g\in G_p}f_1(g)\overline{f_2(g)}.
\]
\end{lemma}

\begin{proof}
For the given class function $f$ on $G_p$,
define the twisted Chebyshev sum by
\[
\Psi_{f,p}(X)
:=
\sum_{\substack{\gamma\in \mathcal P(\Gamma(1)),\ k\geq 1\\
                 N(\gamma)^k\leq X}}
f\!\left(\pi_p(\gamma)^k\right)\log N(\gamma).
\]
By \cite[Theorem~1.2]{Reche2026}, we have
\begin{equation}\label{eq:Reche}
\Psi_{f,p}(X)=  \langle f,\mathbf 1\rangle_{G_p} X + O_{p,f,\epsilon}( X^{25/36+\epsilon}).   
\end{equation}
We put
\[
\Theta_{f,p}(X)
:=
\sum_{\substack{\gamma\in \mathcal P(\Gamma(1)),\\
                 N(\gamma)\leq X}}
f\!\left(\pi_p(\gamma)\right)\log N(\gamma).
\]
Then, we have
\[
\Psi_{f,p}(X)
= \Theta_{f,p}(X) + R_{f,p}(X),
\]
where 
\[
R_{f,p}(X):=\sum_{\substack{\gamma\in \mathcal P(\Gamma(1)),\ k\geq 2\\
                 N(\gamma)^k\leq X}}
f\!\left(\pi_p(\gamma)^k\right)\log N(\gamma).
\]

Put $\ell_0=2\log((3+\sqrt5)/2)$, the minimum primitive
geodesic length on $\Gamma(1)\backslash\mathbb H$.
Since $N(\gamma)^k\le X$ implies $k\ell_0\le\log X$, we have
\[
\begin{aligned}
|R_{f,p}(X)|
&\le \|f\|_\infty
\sum_{2\le k\le(\log X)/\ell_0}
\sum_{\substack{\gamma\in\mathcal P(\Gamma(1))\\
                 N(\gamma)\le X^{1/k}}}
\log N(\gamma)\\
&\ll \|f\|_\infty
\sum_{2\le k\le(\log X)/\ell_0}X^{1/k}\\
&\ll \|f\|_\infty X^{1/2}\log X.
\end{aligned}
\]
Here we used the standard bound
\[
\sum_{\substack{\gamma\in\mathcal P(\Gamma(1))\\
                 N(\gamma)\le Y}}
\log N(\gamma)\ll Y,
\qquad Y\ge2,
\]
which follows from \eqref{eq:Reche} with $f=1$.
Hence,
\begin{align*}\label{eq:Theta-est}
\Psi_{f,p}(X) =\,  \Theta_{f,p} (X) + O_{p,f,\epsilon} (X^{1/2+\epsilon}).
\end{align*}
Combining this with \eqref{eq:Reche}, we have
\begin{equation}\label{eq:est-36-25}
\Theta_{f,p}(X) =   \langle f,\mathbf 1\rangle_{G_p} X + O_{p,f,\epsilon}( X^{25/36+\epsilon}). 
\end{equation}
Now we relate the counting function
\[
\Pi_{f,p}(X)
:=
\sum_{\substack{\gamma\in \mathcal P(\Gamma(1))\\N(\gamma)\leq X}}
f\!\left(\pi_p(\gamma)\right)
\]
with $\Theta_{f,p}(X)$ by
\[
\Pi_{f,p}(X)
=
\frac{\Theta_{f,p}(X)}{\log X}
+
\int_2^X
\frac{\Theta_{f,p}(t)}
     {t(\log t)^2}\,dt.
\]
Substituting \eqref{eq:est-36-25} into this and using the following equality
\[
\operatorname{Li}(X)=\int^X_2 \frac{dt}{\log t} = \frac{X}{\log X} +\int^X_2 \frac{dt}{(\log t)^2} - \frac{2}{\log 2}
\]
yields the stated asymptotic, after absorbing the constant
$2\langle f,1\rangle_{G_p}/\log 2$ into the error term:
\begin{equation*}
\Pi_{f,p}(X)
=
 \langle f,\mathbf 1\rangle_{G_p} \operatorname{Li}(X)
+ O_{p,f,\epsilon}( X^{25/36+\epsilon}).
\label{eq:twisted-prime-geodesic}
\end{equation*}

\end{proof}

We now apply Lemma~\ref{lem:finite-quotient-pgt} to the characteristic
functions of the conjugacy-invariant subsets of \(G_p\) corresponding
to the three noninert reduction types. 
By \eqref{eq:weights-by-reduction-type} and \eqref{eq:Ap-type-decomposition}, we have
\begin{equation}\label{eq:even-reduction-decomp}
A_p^+(x)
=
A_p^{\Split,+}(x)
+
A_p^{\Ramified,+}(x)
+
A_p^{\Scalar,+}(x),
\end{equation}
and
\begin{equation*}\label{eq:odd-reduction-decomp}
A_p^-(x)
=
A_p^{\Ramified,-}(x)
+
A_p^{\Scalar,-}(x).
\end{equation*}
The following theorem determines the leading term of each
even reduction type in \eqref{eq:even-reduction-decomp} and shows that these leading terms cancel exactly.

\begin{theorem}
\label{thm:fixed-level-cancellation}
For every fixed prime $p \ge 5$, as $x\to\infty$,
\[
\sum_{d\in \Dset_{0,p}^{\Split}(x)}h(d)
 \sim
 \frac{p-3}{2(p-1)}\operatorname{Li}(x^2),
\]
\[
\sum_{d\in \Dset_{0,p}^{\Ramified}(x)}h(d)
 \sim
 \frac{2}{p}\operatorname{Li}(x^2),
\]
and
\[
\sum_{d\in \Dset_{0,p}^{\Scalar}(x)}h(d)
 \sim
 \frac{2}{p(p^2-1)}\operatorname{Li}(x^2).
\]
Consequently,
\[
A_p^{\Split,+}(x)
 \sim
 -\frac{p-3}{p-1}\operatorname{Li}(x^2),
\]
\[
A_p^{\Ramified,+}(x)
 \sim
 \frac{p-3}{p}\operatorname{Li}(x^2),
\]
and
\[
A_p^{\Scalar,+}(x)
 \sim
 \frac{p-3}{p(p-1)}\operatorname{Li}(x^2).
\]
The three leading coefficients cancel exactly.
Hence, by \eqref{eq:even-reduction-decomp},
\[
\begin{aligned}
A_p^+(x)
 &=
 A_p^{\Split,+}(x)
 +A_p^{\Ramified,+}(x)
 +A_p^{\Scalar,+}(x) \\
 &=
 o\!\left(\operatorname{Li}(x^2)\right).
\end{aligned}
\]
\end{theorem}

\begin{proof}
Let
\[
T_{\Split},\quad
T_{\Ramified},\quad
T_{\Scalar} \qquad \subset \quad G_p
\]
denote the subsets corresponding to the split, ramified, and scalar reduction types, respectively.
The number of split maximal tori in $G_p$ is
\[
[G_p: N_{G_p}(T)]=\frac{|G_p|}{p-1}=\frac{p(p^2-1)}{2} \cdot \frac{1}{p-1}  =\frac{p(p+1)}2,
\]
where $N_{G_p}(T)$ is the normalizer of a fixed torus $T$, whose order is $p-1$.
Each such torus contains \((p-3)/2\) nonidentity elements, and every element of split reduction type belongs to a unique split maximal torus.
Hence
\[
|T_{\Split}|
 =\frac{p(p+1)(p-3)}4.
\]
Similarly, \(G_p\) has \(p+1\) Sylow \(p\)-subgroups, each containing \(p-1\) nonidentity elements. These elements are precisely those of ramified reduction type, and each belongs to a unique Sylow \(p\)-subgroup. Therefore
\[
|T_{\mathrm{ram}}|=(p+1)(p-1)=p^2-1.
\]
The scalar reduction type corresponds to the identity of \(G_p\), so \(|T_{\mathrm{sc}}|=1\).
Their densities in $G_p$ are consequently
\[
q_{\Split}
 =\frac{|T_{\Split}|}{|G_p|}
 =\frac{p-3}{2(p-1)},
\qquad
q_{\Ramified}
 =\frac{|T_{\Ramified}|}{|G_p|}
 =\frac{2}{p},
\qquad
q_{\Scalar}
 =\frac{|T_{\Scalar}|}{|G_p|}
 =\frac{2}{p(p^2-1)}.
\]

Let $I_T$ denote the characteristic function of $T$, where $T$ is one of the three subsets
$T_{\Split}$, $T_{\Ramified}$, and $T_{\Scalar}$ of $G_p$.
By Lemma~\ref{lem:finite-quotient-pgt},
\[
\sum_{\substack{\gamma\in \mathcal P(\Gamma(1))\\N(\gamma)\leq X}}
 I_T\bigl(\pi_p(\gamma)\bigr)
 =
 q_T\operatorname{Li}(X)
 +o\!\left(\operatorname{Li}(X)\right).
\]
For a discriminant $d\in \Dset$, the corresponding primitive
geodesics have norm $N_d=\epsilon_d^2$,
and their class number is $h(d)$.  Thus the condition $\epsilon_d\leq x$ is
equivalent to $N_d\leq x^2$.  By the reduction type classification of Section~\ref{sec:class-number-geodesics},
the images in \(G_p=\mathrm{PSL}_2(\mathbb F_p)\) of the split, ramified, and scalar discriminants belong to \(T_{\Split}\), \(T_{\Ramified}\), and \(T_{\Scalar}\), respectively.
Hence
\[
\sum_{d\in \Dset_{0,p}^{\Split}(x)}h(d)
 =
 \sum_{\substack{\gamma\in \mathcal P(\Gamma(1))\\N(\gamma)\leq x^2}}
 I_{T_{\Split}}\bigl(\pi_p(\gamma)\bigr),
\]
and similarly for the ramified and scalar types.  Applying the preceding
asymptotic formulas with $X=x^2$ proves
\[
\sum_{d\in \Dset_{0,p}^{\Split}(x)}h(d)
 \sim
 \frac{p-3}{2(p-1)}\operatorname{Li}(x^2),
\]
\[
\sum_{d\in \Dset_{0,p}^{\Ramified}(x)}h(d)
 \sim
 \frac{2}{p}\operatorname{Li}(x^2),
\]
and
\[
\sum_{d\in \Dset_{0,p}^{\Scalar}(x)}h(d)
 \sim
 \frac{2}{p(p^2-1)}\operatorname{Li}(x^2).
\]

We now apply the nontrivial even-character weights in \eqref{eq:weights-by-reduction-type} to the class-number sums \eqref{eq:Ap-type-components} over the split, ramified, and scalar discriminants.
These weights are, respectively,
\[
-2,\qquad
\frac{p-3}{2},\qquad
\frac{(p+1)(p-3)}{2}.
\]
It follows that
\[
A_p^{\Split,+}(x)
 \sim
 -\frac{p-3}{p-1}\operatorname{Li}(x^2),
\]
\[
A_p^{\Ramified,+}(x)
 \sim
 \frac{p-3}{p}\operatorname{Li}(x^2),
\]
and
\[
A_p^{\Scalar,+}(x)
 \sim
 \frac{p-3}{p(p-1)}\operatorname{Li}(x^2).
\]

Their leading coefficients satisfy
\[
-\frac{p-3}{p-1}
+\frac{p-3}{p}
+\frac{p-3}{p(p-1)}
=0,
\]
and thus we have
\[
A_p^+(x)
 =o\!\left(\operatorname{Li}(x^2)\right).
\]
\end{proof}

\begin{remark}
The \(\sigma_p(d)\)-weighted scalar and ramified contributions occur in the odd-character sum \(A_p^-(x)\) and do not descend to functions on the projective quotient \(G_p\). Consequently, the argument used in Lemma~\ref{lem:finite-quotient-pgt} cannot be applied separately to these contributions.
An analogous conclusion for \(A_p^-(x)\) would
require signed equidistribution between the classes with
\(t_d\equiv2\pmod p\) and \(t_d\equiv-2\pmod p\), separately within the
scalar and ramified types. 
\end{remark}

\section{Numerical behavior of the reduction type contributions}

\label{sec:local-type-contributions}

Section~\ref{sec:class-number-geodesics} decomposed the character sums into scalar, split, ramified,
and inert reduction types, while the experiments of Section~\ref{sec:experimental-setup}
suggested cancellation on the scale \(x\sqrt p\). Section~\ref{sec:cancellation-thm} then
determined the \(\operatorname{Li}(x^2)\) leading term of each even component and proved that these leading terms cancel exactly. This result, however, does not describe the remaining
level-dependent fluctuations. 

In this section, we first illustrate in Figure~\ref{fig:split-ramified} the striking symmetry between the split and ramified contributions. We then subtract the predicted leading terms from the scalar, split, and ramified contributions and study the resulting errors numerically. The purpose of this section is to determine which reduction types account for the observed \(x\sqrt p\)-scale and to examine how the symmetry and correlated fluctuations among these contributions produce the overall oscillatory pattern.

\subsection{Normalization and binning}\label{subsection:normalization}

Theorem~\ref{thm:fixed-level-cancellation} shows that the individual split,
ramified, and scalar sums contain leading terms of order
\(\operatorname{Li}(x^2)\) respectively, although these terms cancel in their even
combination. 
After normalization by \(x\sqrt p\), an individual leading term has size
\[
\frac{\operatorname{Li}(x^2)}{x\sqrt p}
\asymp
\frac{x}{\sqrt p\log x}.
\]
The purpose of the first experiment in this section is to illustrate
the cancellation between these contributions.
We then subtract the main term of each reduction type component and
examine the remaining fluctuations using the normalization \(x\sqrt p\).

To study the simultaneous dependence on the level $p$ and the length $x$, we organize the data in terms of the parameter
\[
\eta=\frac{p}{x}.
\]
Restricting \(\eta\) to a short interval compares pairs \((p,x)\) in a
similar relative level range while retaining sufficiently many primes
for stable bin averages. The interval specified below is used solely as a sampling
window and has no intrinsic significance.
The main sample window used in the present figures is
\[
  E_0=[0.002,0.006].
\]
The sampled \(x\)-values are
\[
  X_{\rm samp}=
  \{\,x\in \Z:\ 200000\leq x\leq10^6,\ x\equiv0\pmod {5000},\
  P_x(E_0)\neq\emptyset\,\},
\]
where
\[
  P_x(E_0)=\{\,p\in \Z_{>0}:\ p\text{ odd prime},\  \eta \in E_0\,\}.
\]
This gives $48,963$ individual $(p,x)$-pairs, with $|X_{\rm samp}|=161$.

For the bin-mean curves in Figure~\ref{fig:split-ramified} and Figures~\ref{fig:survivor-split}--\ref{fig:survivor-scalar-sigma}, we divide $E_0$ into $24$ equal-width intervals.
In each interval, we plot the arithmetic mean of the corresponding $y$-values at the midpoint of the interval, joining these points by line segments.

\subsection{Cancellation between the split and ramified contributions}

\begin{figure}[htbp]
\centering
\includegraphics[width=0.78\textwidth]{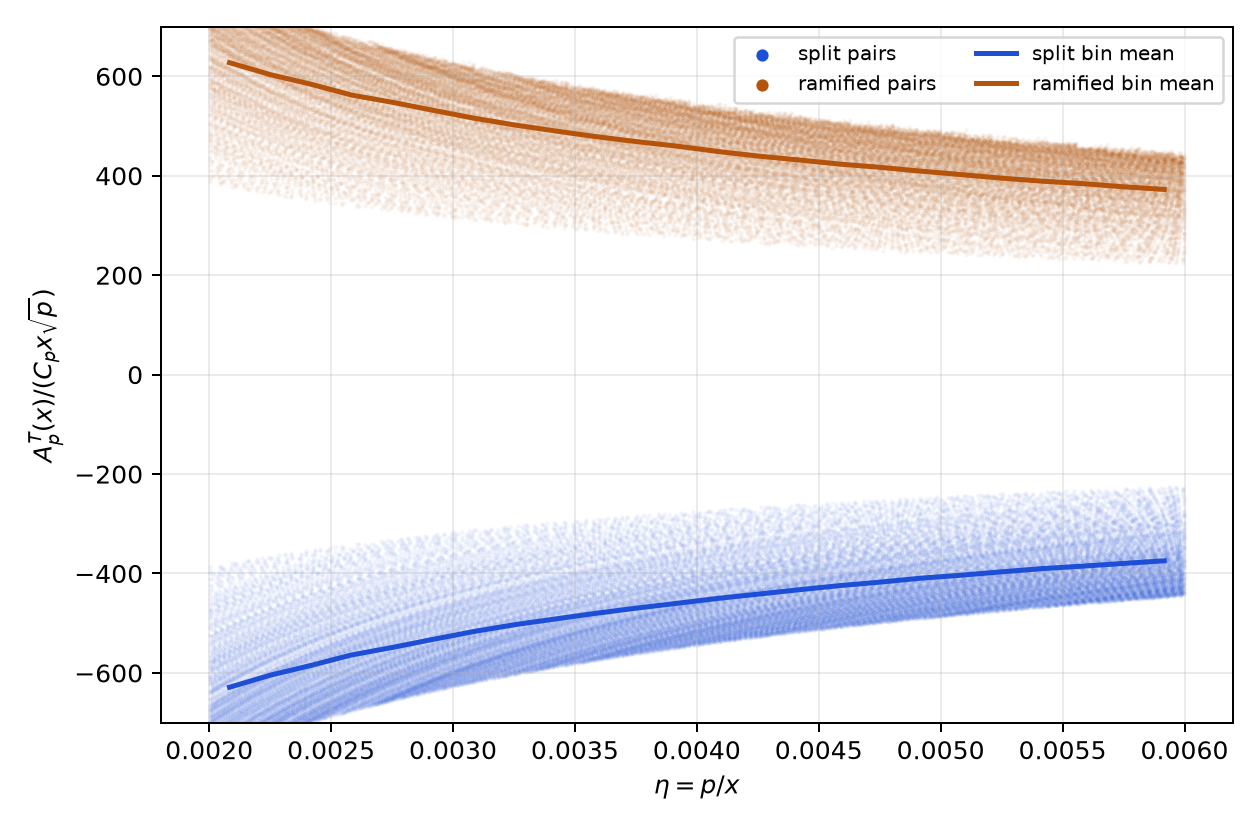}
\caption{\sf \small \(x\)-axis: \(\eta=p/x\).  \(y\)-axis: normalized local type sums.
The split contribution is large and negative, while the ramified contribution is large and positive.
The total \(A_p(x)\) is a comparatively small remainder after this cancellation, together with the
scalar contribution.}
\label{fig:split-ramified}
\end{figure}

Figure~\ref{fig:split-ramified}  compares the split and ramified contributions and illustrates the cancellation between their leading terms. The split contribution is entirely contained in the even character part, whereas the ramified contribution decomposes as

$$ A_p^{\Ramified,+}(x) + A_p^{\Ramified,-}(x). $$

The ramified curve in Figure~\ref{fig:split-ramified} is dominated by the positive \(\Li(x^2)\)-main term of \(A_p^{\Ramified,+}(x)\). As Figure~\ref{fig:survivor-ramified-sigma} shows, the odd ramified contribution fluctuates around zero on the \(x\sqrt p\)-scale and has no comparable leading profile. It is therefore negligible on the much larger scale used in Figure~\ref{fig:split-ramified}. Thus Figure~\ref{fig:split-ramified} primarily displays the cancellation between the leading even part terms of the split and ramified contributions. The smaller scalar leading term completes the exact cancellation established in Theorem~\ref{thm:fixed-level-cancellation}.

\subsection{Reduction type error terms}
Define the five errors obtained by subtracting the leading terms:
\begin{equation}
\label{eq:local-type-main-subtracted-quantities}
\begin{aligned}
  B_p^{\Split,+}(x)
  &=
  \frac{
    A_p^{\Split,+}(x)
    +\frac{p-3}{p-1}\Li(x^2)
  }{x\sqrt p},\\
  B_p^{\Ramified,+}(x)
  &=
  \frac{
    A_p^{\Ramified,+}(x)
    -\frac{p-3}{p}\Li(x^2)
  }{x\sqrt p},\\
  B_p^{\Scalar,+}(x)
  &=
  \frac{
    A_p^{\Scalar,+}(x)
    -\frac{p-3}{p(p-1)}\Li(x^2)
  }{x\sqrt p},\\
  B_p^{\Ramified,-}(x)
  &=
  \frac{A_p^{\Ramified,-}(x)}{x\sqrt p},
  \qquad
  B_p^{\Scalar,-}(x)
  =
  \frac{A_p^{\Scalar,-}(x)}{x\sqrt p}.
\end{aligned}
\end{equation}
Because the three main terms cancel exactly, we have
\[
B_p^{\Split,+}(x)
+B_p^{\Ramified,+}(x)
+B_p^{\Scalar,+}(x)
=
\frac{A_p^+(x)}{x\sqrt p}.
\]
Similarly,
\[
B_p^{\Ramified,-}(x)
+B_p^{\Scalar,-}(x)
=
\frac{A_p^-(x)}{x\sqrt p}.
\]

Thus bounds for the errors associated with individual reduction
types yield corresponding bounds for the even, odd, and combined
character sums. Each quantity in
\eqref{eq:local-type-main-subtracted-quantities} measures the
deviation from its leading term, normalized by \(x\sqrt p\).

\begin{figure}[htbp]
\centering
\includegraphics[width=0.78\textwidth]{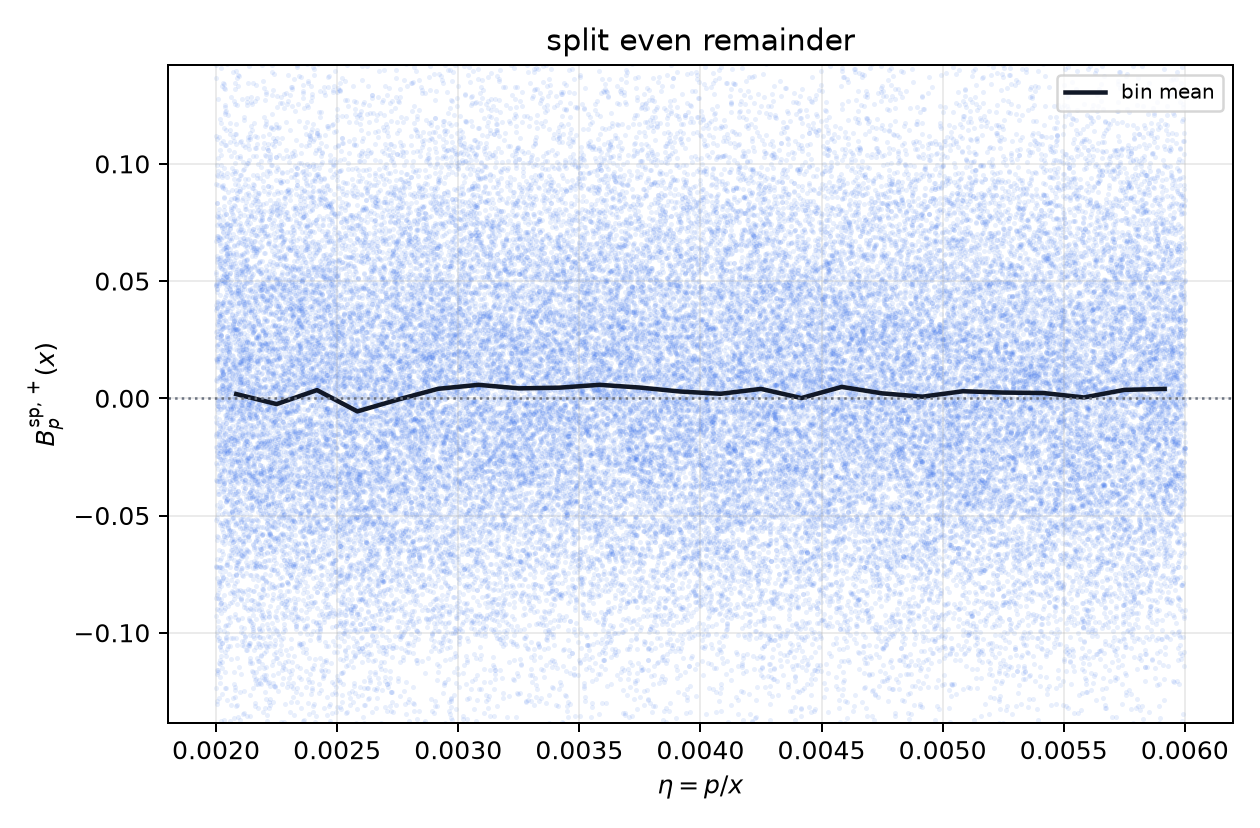}
\caption{\sf \small \(B_p^{\Split,+}(x)\) against \(\eta=p/x\). The solid curve shows the bin means.}
\label{fig:survivor-split}
\end{figure}

\begin{figure}[htbp]
\centering
\includegraphics[width=0.78\textwidth]{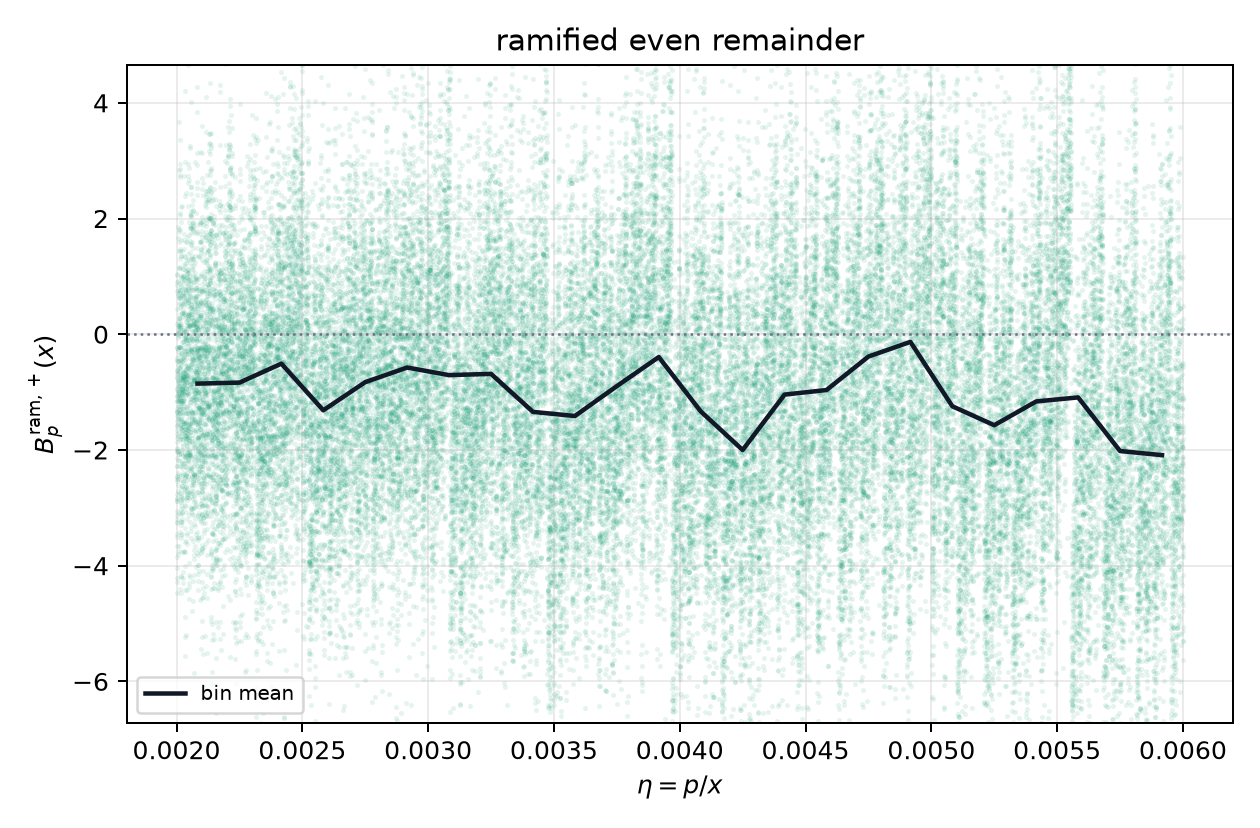}
\caption{\sf \small \(B_p^{\Ramified,+}(x)\) against \(\eta=p/x\). The solid curve shows the bin means.}
\label{fig:survivor-ramified-zero}
\end{figure}

\begin{figure}[htbp]
\centering
\includegraphics[width=0.78\textwidth]{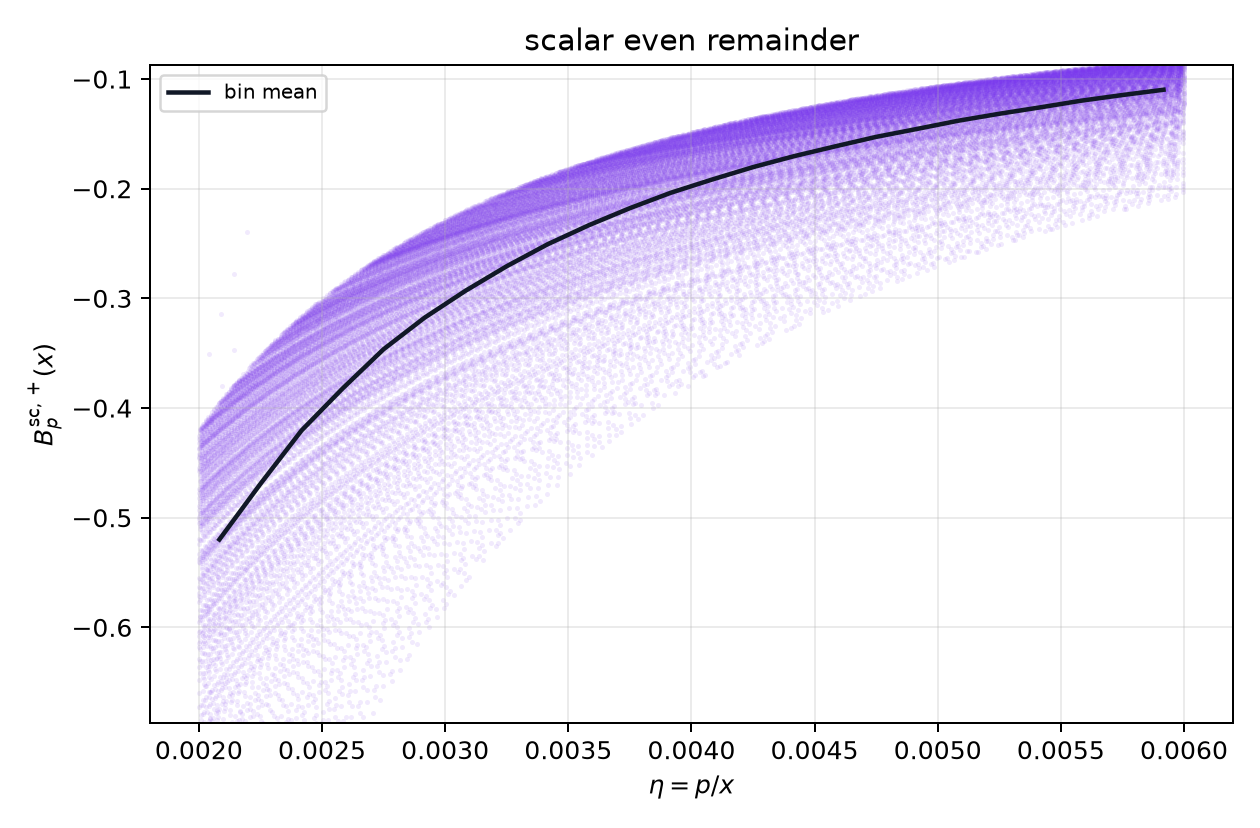}
\caption{\sf \small \(B_p^{\Scalar,+}(x)\) against \(\eta=p/x\). The solid curve shows the bin means.}
\label{fig:survivor-scalar-zero}
\end{figure}

\begin{figure}[htbp]
\centering
\includegraphics[width=0.78\textwidth]{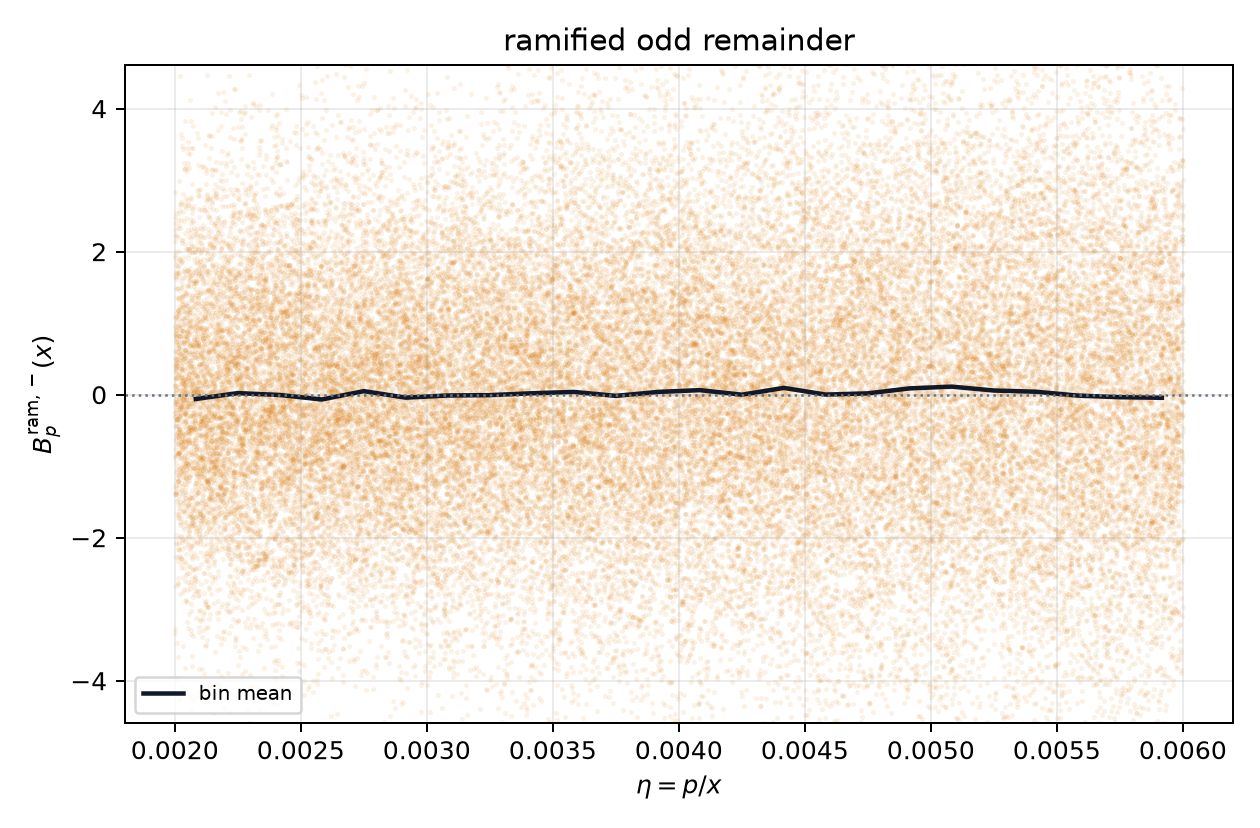}
\caption{\sf \small \(B_p^{\Ramified,-}(x)\) against \(\eta=p/x\).  This plot measures the
odd-sign ramified fluctuation after normalization by \(x\sqrt p\).}
\label{fig:survivor-ramified-sigma}
\end{figure}

\begin{figure}[htbp]
\centering
\includegraphics[width=0.78\textwidth]{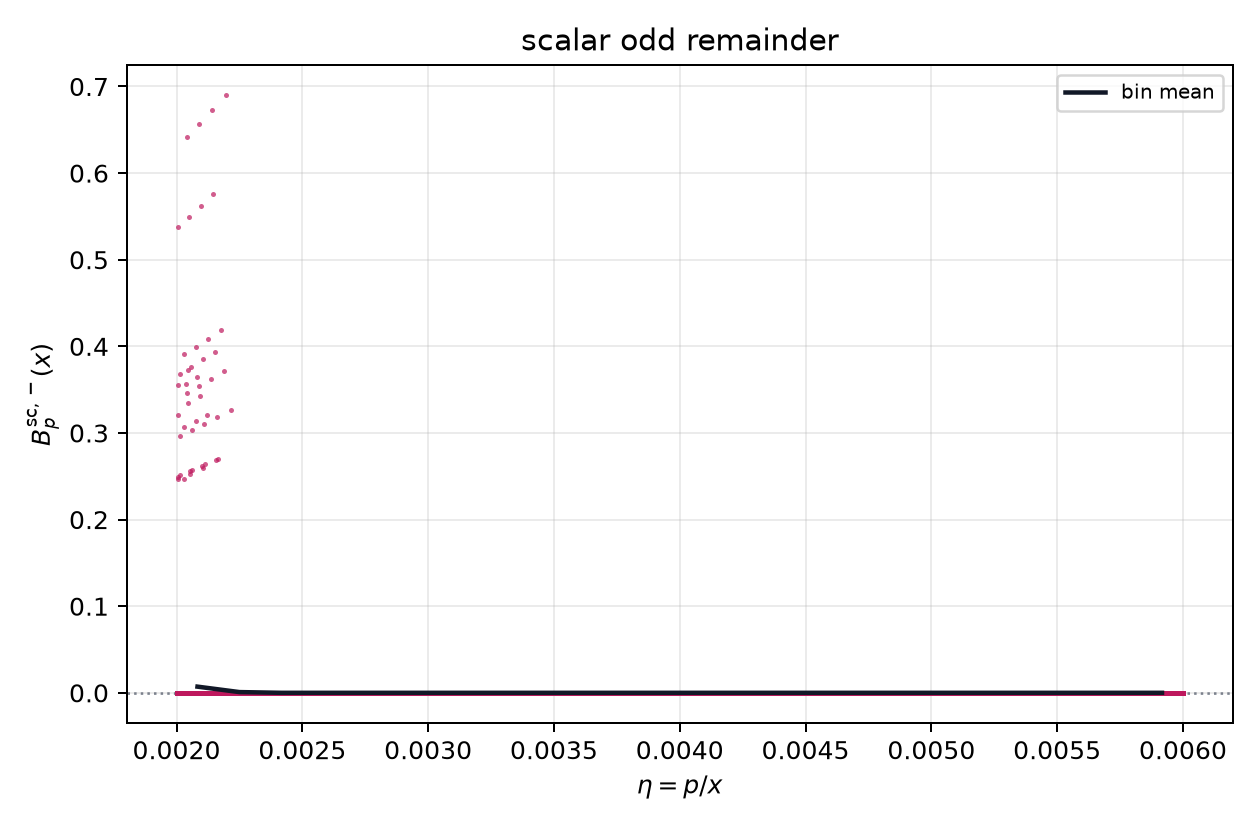}
\caption{\sf \small \(B_p^{\Scalar,-}(x)\) against \(\eta=p/x\).  This is the scalar
odd-sign term after normalization; it is small in comparison with the raw scalar coefficient.}
\label{fig:survivor-scalar-sigma}
\end{figure}

\FloatBarrier

Figures~\ref{fig:survivor-split}--\ref{fig:survivor-scalar-sigma}
show that the five contributions, after subtraction of their
leading terms, exhibit distinct numerical behavior but remain
bounded on the $x\sqrt p$ scale throughout the
sampled fixed-ratio window.
The split
remainder in Figure~\ref{fig:survivor-split} is concentrated near zero, while the ramified
$0$-part in Figure~\ref{fig:survivor-ramified-zero} exhibits fluctuations whose bin means remain
of bounded size.  The scalar $0$-part in Figure~\ref{fig:survivor-scalar-zero} displays a coherent
negative profile, discussed separately below in Remark \ref{rmk:interval}. The ramified $\sigma_p$-part in Figure~\ref{fig:survivor-ramified-sigma} fluctuates
around zero, whereas the scalar $\sigma_p$-part in Figure~\ref{fig:survivor-scalar-sigma} is
comparatively small.
Taken together, these observations suggest bounded behavior of the normalized remainders across the sampled primes and cutoffs in the fixed-ratio window. The arithmetic jumps discussed above prevent extending this interpretation to a uniform bound at every sharp cutoff.

\begin{remark}\label{rmk:interval}
The patterns in Figures \ref{fig:survivor-scalar-zero} and \ref{fig:survivor-scalar-sigma} have a common origin in the chosen sampling window. In the subsection \ref{subsection:normalization}, it is assumed that \(p\ge0.002x\) and \(x\ge200000\), hence \(x/p^2\le1.25\). Scalar reduction requires \(t_d\equiv\pm2\pmod{p^2}\). The candidate \(p^2-2\) comes from squaring a unit of trace \(p\), so the only primitive scalar trace accessible here is \(p^2+2\), with positive sign.
Put \(t_0=p^2+2\). Within the sampled window, every scalar discriminant contributing to the sum has \(t_d=t_0\), and hence
\[
\epsilon_d=\alpha_{t_0}
=\frac{t_0+\sqrt{t_0^2-4}}2.
\]
Thus all these discriminants enter the sum with the condition $x \ge \alpha_{t_0}$. 

Define their class-number sum by
\[
M_p=\sum_{\substack{d\in\mathcal D:\,t_d=t_0\\p\mid u_d}}h(d).
\]
Since \(p\mid u_d\) and \(t_d\equiv2\pmod p\), the corresponding matrix reduces to \(I\). Each induced representation associated with an odd character therefore contributes its dimension \(p+1\). Summing over the \((p-1)/2\) odd characters gives
\[
A_p^{\mathrm{sc},-}(x)=
\begin{cases}
0,&x<\alpha_{t_0},\\[2pt]
\dfrac{p^2-1}{2}M_p,&x\ge\alpha_{t_0},
\end{cases}
\]
for sampled pairs \((p,x)\). Consequently, when \(x\ge\alpha_{t_0}\),
\[
B_p^{\mathrm{sc},-}(x)
=\frac{A_p^{\mathrm{sc},-}(x)}{x\sqrt p}
=\frac{(p^2-1)M_p}{2p^{3/2}}\,\eta,
\qquad \eta=\frac px.
\]
The coefficient of \(\eta\) depends only on \(p\), explaining the straight rows of positive dots in Figure~\ref{fig:survivor-scalar-sigma}.
Only \(47\) of the \(48\,963\) sampled pairs satisfy the condition $x \ge \alpha_{t_0}$ and therefore have a nonzero scalar contribution. 

The even scalar weight is \((p+1)(p-3)/2\). Thus, throughout this window,
\(A_p^{\mathrm{sc},+}=\frac{p-3}{p-1}A_p^{\mathrm{sc},-}\).
Applying the definitions in \eqref{eq:local-type-main-subtracted-quantities} gives
\[
B_p^{\mathrm{sc},+}(x)
=
\frac{p-3}{p-1}B_p^{\mathrm{sc},-}(x)
-\frac{p-3}{p(p-1)}
\frac{\operatorname{Li}(x^2)}{x\sqrt p}.
\]
For sampled pairs satisfying \(x<\alpha_{p^2+2}\), the scalar sums vanish, so the even remainder consists solely of the negative subtracted main term.
Since \(\operatorname{Li}(x^2)\sim x^2/(2\log x)\),
\[
B_p^{\mathrm{sc},+}(x)
\approx
-\frac{p-3}{p-1}
\frac{\eta^{-3/2}}{2\sqrt{x}\log x}.
\]
This accounts for the negative profile in Figure~\ref{fig:survivor-scalar-zero}. 
The nonzero scalar contributions affect only 47 sampled pairs and have little visible effect on Figure~\ref{fig:survivor-scalar-zero}, whereas they appear as positive points in Figure~\ref{fig:survivor-scalar-sigma}.
\end{remark}

\section{\texorpdfstring
  {From $A^+_p(x)$-estimate to a spectral gap}
  {From A-plus(x)-estimate to a spectral gap}}
\label{sec:dirichlet-spectral-gap}

In this section, we derive the conditional spectral consequence of $A^+_p(x)$-estimate in Conjecture~\ref{conj:ap-experimental}. Weakening the inequality in Conjecture~\ref{conj:ap-experimental}, we consider the following hypothesis:

\begin{hypothesis}
\label{hyp:trace-A+p}
For $B >0$ and $\epsilon >0$,
there exist constants $C=C(B)$ and $\nu=\nu(B) \geq 0$ such that
\begin{equation}
|A_p^+(x)|
\ll_{B,\epsilon}
p^C x^{1+\nu} (p x)^{\epsilon}
\qquad \text{for} \quad
  x\leq p^B, 
\label{eq:arithmetic-input}
\end{equation}
uniformly in $p$.
\end{hypothesis}

The passage from \eqref{eq:arithmetic-input} to a spectral counting estimate proceeds in three steps. First, we use \eqref{eq:arithmetic-input} to bound the primitive hyperbolic contribution to the trace formula. Second, we choose a test function whose spectral transform is nonnegative and amplifies the contribution of exceptional eigenvalues. Third, we bound the remaining terms in the trace formula. Combining these estimates yields a bound for the number of exceptional eigenvalues.

\subsection{The primitive hyperbolic contribution}
To incorporate the primitive orbital weight, recall that
\begin{equation}
\label{eq:A-p^+}
A_p^+(x)
=
\sum_{d\in  \Dset_{0,p}(x)} a_p^+(d),
\end{equation}
and define its weighted analogue by
\begin{equation*}
 A_{p,F}^+(x)
 :=
 \sum_{d\in  \Dset_{0,p}(x)} a_p^+(d)F(\epsilon_d),
 \qquad
 F(u):=\frac{2\log u}{1-u^{-2}}.
 \label{eq:weighted-even-sum}
\end{equation*}
By the definition in \eqref{eq:A-p^+},
the weighted sum is the Stieltjes integral
\[
 A_{p,F}^+(x)=\int_{1^-}^{x}F(u)\,dA_p^+(u).
\]
Integration by parts therefore gives
\[
 A_{p,F}^+(x)
 =
 A_p^+(x)F(x)
 -
 \int_1^x A_p^+(u)F'(u)\,du.
\]
Since
\[
F(u)\ll 1+ \log u,
\qquad
F'(u)\ll u^{-1} \qquad \text{for} \quad u\geq 1,
\]
the arithmetic assumption \eqref{eq:arithmetic-input} gives
\begin{equation}
A_{p,F}^+(x)
\ll_{B,\epsilon}\,
p^C x^{1+\nu}\, (\log x)\, (px)^{\epsilon}
\qquad
\text{for} \quad x\leq p^B.
\label{eq:weighted-arithmetic-bound}
\end{equation}
We next choose a positive family of test functions.

\begin{lemma}\label{lem:positive-test-family}
Fix $0<\kappa<1$. There exist even, nonnegative functions
$g_T\in C_c^\infty(\mathbb R)$, $T\ge1$, whose Fourier transforms
$h_T(z):=\int_{\mathbb R}g_T(u)e^{izu}\,du$ satisfy:
\begin{enumerate}
\item $\operatorname{supp}(g_T)\subset[-T,T]$ and, for every $j\ge0$,
\[
 \|g_T^{(j)}\|_\infty\ll_{j,\kappa}T^{-j},
 \qquad \|g_T\|_{1}\ll_\kappa T.
\]
\item For every $N\ge0$ and $r\in\mathbb R$,
\begin{equation} \label{eq:test-seminorms-remainder}
 0\le h_T(r)\ll_{N,\kappa}T(1+T|r|)^{-N}.
\end{equation}
Consequently, $h_T(0)\ll_\kappa T$ and, for every fixed $a\ge0$,
\begin{equation} \label{eq:hT-moments-remainder}
 \int_{\mathbb R}(1+|r|)^a h_T(r)\,dr\ll_{a,\kappa} 1.
\end{equation}
\item Uniformly for $0<\eta<\tfrac12$,
\[
 h_T(i\eta)\gg_\kappa T e^{\kappa\eta T}.
\]
In particular, given $B>0$, $0<r_0<\tfrac12$, and $\epsilon>0$,
one may fix $\kappa<1$ sufficiently close to $1$ so that
\[
 h_{2B\log p}(i\eta)\gg_{B,r_0,\epsilon}
 p^{2Br_0-\epsilon}
 \qquad(r_0\le\eta<\tfrac12,\;2B\log p\ge1).
\]
\end{enumerate}
\end{lemma}

\begin{proof}
Choose an even, nonnegative
$\phi\in C_c^\infty((-1/2,1/2))$ with $\phi(\kappa/2)>0$,
and set $G=\phi*\phi$ and $g_T(u)=G(u/T)$. Equivalently,
$g_T=T^{-1}(\phi_T*\phi_T)$, where $\phi_T(u)=\phi(u/T)$.
The support, nonnegativity, derivative, and $L^1$ bounds
follow directly by scaling. Moreover,
\[
h_T(z)=T\widehat\phi(Tz)^2.
\]
Since $\widehat\phi(r)$ is real for real $r$,
we have $h_T(r)\ge0$, and its rapid decay gives
\eqref{eq:test-seminorms-remainder}.
Substituting $v=Tr$ and taking $N>a+1$ yields
\eqref{eq:hT-moments-remainder}, since $T\ge1$.

For the third claim, choose $\delta>0$ such that
$c_\kappa:=\int_{\kappa/2}^{\kappa/2+\delta}\phi(v)\,dv>0$.
By evenness and nonnegativity,
\[
h_T(i\eta)
=T\left(\int_{\mathbb R}\phi(v)e^{\eta Tv}\,dv\right)^2
\ge c_\kappa^2 T e^{\kappa\eta T}.
\]
Finally, given $B,r_0,\epsilon$, choose $\kappa<1$
so that $2Br_0(1-\kappa)<\epsilon$.
Substituting $T=2B\log p\ge1$ and using $\eta\ge r_0$
gives the desired bound.
\end{proof}

\subsection{The summed trace formula}

We sum the \(\chi_{\Gamma_0(p)}\)-twisted Selberg trace formulas on \(\Gamma_0(p)\) over all nontrivial even characters modulo \(p\). By induction, this sum is equivalent to the vector-valued Selberg trace formula on \(\Gamma(1)\) associated with
\begin{equation}\label{eq:induced-rep-Gamma}
\rho_{\mathrm{Dir},\Gamma}
:=\rho_{\mathrm{Dir}}\circ\pi_p
\cong
\bigoplus_{\substack{\chi\bmod p\\ \chi(-1)=1,\ \chi\ne \mathbf 1}}
\operatorname{Ind}_{\Gamma_0(p)}^{\Gamma(1)}
\chi_{\Gamma_0(p)},
\end{equation}
where \(\rho_{\mathrm{Dir}}\) is defined in \eqref{eq:def-rho-Dir}  and \(\pi_p:\Gamma(1)\to G_p\) is the reduction map. 
The induction correspondence gives a unitary identification
\[
\bigoplus_{\substack{\chi\bmod p\\
\chi(-1)=1,\ \chi\ne \mathbf 1}}
L^2_{\mathrm{cusp}}
\bigl(\Gamma_0(p)\backslash\mathbb H,\chi_{\Gamma_0(p)}\bigr)
\cong
L^2_{\mathrm{cusp}}
\bigl(\Gamma(1)\backslash\mathbb H,\rho_{\mathrm{Dir},\Gamma}\bigr),
\]
which intertwines the Laplacians and preserves spectral multiplicities. Hence the exceptional eigenvalues counted by \(N_{\mathrm{Dir}}(r;p)\) are precisely those of the vector-valued Laplacian associated with \((\Gamma(1),\rho_{\mathrm{Dir},\Gamma})\). Although the base surface is now \(Y(1)\), the level and character dependence are retained in the representation \(\rho_{\mathrm{Dir},\Gamma}\).

Put \(X=e^{T/2}\). We restrict \(1\le T\le2B\log p\), so that \(X\le p^B\) and the arithmetic estimate applies throughout the range of integration.
Let \(r_{\chi,j}\) denote the spectral parameters of the discrete Laplace spectrum on \(L^2(\Gamma_0(p)\backslash\mathbb H,\chi_{\Gamma_0(p)})\), counted with multiplicity. By the induction correspondence underlying \eqref{eq:induced-rep-Gamma}, the direct sum of these spectra over all nontrivial even characters is the discrete spectrum of the vector-valued Laplacian associated with \(\rho_{\mathrm{Dir},\Gamma}\). Applying the test pair \((g_T,h_T)\) of Lemma \ref{lem:positive-test-family} and separating the primitive hyperbolic contribution in this vector-valued trace formula,
we obtain
\begin{equation}
\sum_{\substack{\chi\bmod p\\
                \chi(-1)=1,\ \chi\neq \mathbf 1}}
\ \sum_j h_T(r_{\chi,j})
=
H_{\mathrm{Dir},1}(T)+R_{\mathrm{Dir}}(T).
\label{eq:summed-trace-formula}
\end{equation}
The spectral sum on the left hand side of  \eqref{eq:summed-trace-formula} runs over the full discrete spectrum,
including any residual eigenvalues. Since $h_T$ is nonnegative
on both the real and exceptional imaginary spectral parameters,
this sum dominates the contribution of the cuspidal eigenvalues
counted by $N_{\mathrm{Dir}}(r;p)$.
On the right hand side of \eqref{eq:summed-trace-formula},  $R_{\mathrm{Dir}}(T)$ contains all remaining terms, with the
appropriate signs after they have been moved to the right-hand side:
the imprimitive hyperbolic terms, the identity, parabolic and elliptic
terms, and the continuous-spectrum contribution.

Note that the primitive hyperbolic contribution is given  by
\[
 H_{\mathrm{Dir},1}(T)
 =
\sum_{\gamma \in\mathcal P(\Gamma(1))}
 \operatorname{tr}\rho_{\mathrm{Dir},\Gamma}(\gamma)\,
 \frac{\ell(\gamma)}{2\sinh(\ell(\gamma)/2)}\,g_T(\ell(\gamma)).
\]
Here $\mathcal P(\Gamma(1))$ is the set of primitive hyperbolic
conjugacy classes in $\Gamma(1)/\{\pm I\}$, and the subscript
$1$ denotes the term $m=1$ in the sum over powers $\gamma^m$.
Only classes with $\ell(\gamma)\le T$ contribute.
Grouping these classes by discriminant, using
$\ell(\gamma)=2\log\epsilon_d$ and the character-summed
multiplicity $a_p^+(d)$, gives
\begin{equation}\label{eq:primitive-hyperbolic-stieltjes}
 \begin{aligned}
 H_{\mathrm{Dir},1}(T)
 &=
 \sum_{d\in  \Dset_{0,p}(e^{T/2})}
 a_p^+(d)\,
 \frac{2\log\epsilon_d}{\epsilon_d-\epsilon_d^{-1}}\,
 g_T(2\log\epsilon_d)\\
 &=
 \int_{1^-}^{e^{T/2}}
 \frac{g_T(2\log u)}{u}\,dA_{p,F}^+(u),
 \end{aligned}
\end{equation}
since
\[
 \frac{2\log u}{u-u^{-1}}=\frac{F(u)}{u}.
\]
Integrating by parts in
\eqref{eq:primitive-hyperbolic-stieltjes}, with vanishing endpoint terms,
gives
\[
H_{\mathrm{Dir},1}(T)
=
-
\int_1^{e^{T/2}}
A_{p,F}^+(u)\,
d\left(\frac{g_T(2\log u)}{u}\right).
\]
Since
\[
\frac{d}{du}
\left(
\frac{g_T(2\log u)}{u}
\right)
=
\frac{2g_T'(2\log u)-g_T(2\log u)}{u^2},
\]
Lemma~\ref{lem:positive-test-family} and
\eqref{eq:weighted-arithmetic-bound} imply
\begin{equation}\label{eq:primitive-hyperbolic-bound}
\begin{aligned}
|H_{\mathrm{Dir},1}(T)|
&\ll_{B,\varepsilon,\kappa}
p^{C+\varepsilon}
\int_1^X u^{-1+\nu+\varepsilon}\log u\,du\\
&\ll_{B,\varepsilon,\kappa}
p^{C+\varepsilon}X^{\nu+\varepsilon}(1+T)^2,
\qquad X=e^{T/2}.
\end{aligned}
\end{equation}
For the other term $R_{\mathrm{Dir}}(T)$ on the right hand side of \eqref{eq:summed-trace-formula}, we have

\begin{proposition}
\label{hyp:trace-remainder}
For every $\epsilon>0$, 
\begin{equation}\label{eq:R-dir}
|R_{\mathrm{Dir}}(T)|
 \ll_{B,\epsilon,\kappa} p^{2+\epsilon}
\qquad \text{for} \quad
1\leq T\leq 2B \log p.
\end{equation}
\end{proposition}

The proof of this proposition will be given in the next subsection. We can now state the conditional spectral consequence of Hypothesis \ref{hyp:trace-A+p}. 

\begin{proposition}
\label{prop:conditional-dirichlet-gap}
Assume the arithmetic bound given in \eqref{eq:arithmetic-input}.
Then
\begin{equation}\label{eq:densitycount}
 N_{\mathrm{Dir}}(r;p)\ll_{B,\epsilon, r}
 p^{C+B\nu-2Br+\epsilon}+p^{2-2Br+\epsilon}.
\end{equation}
In particular $N_{\mathrm{Dir}}(r;p)=0$ for all sufficiently large
primes if
\begin{equation*}\label{eq:densitythreshold}
 r>\max\left\{\frac\nu2+\frac{C}{2B},\frac1B\right\}.
\end{equation*}
\end{proposition}

\begin{proof}
Take $T=2B\,\log p$, so that $X=p^B$. 
Then, under the assumption in \eqref{eq:arithmetic-input} and
the bounds in \eqref{eq:primitive-hyperbolic-bound} and
\eqref{eq:R-dir}, the right-hand side of
\eqref{eq:summed-trace-formula} is bounded for $x\leq p^B$
as follows:
\[
| H_{\mathrm{Dir},1}(T)| \, +\, | R_{\mathrm{Dir}}(T)| \, \ll_{B,\epsilon,r}\,  p^{C+B\nu+\epsilon/2}+p^{2+\epsilon/2}.
\]
On the other hand, the left hand side of \eqref{eq:summed-trace-formula} has the following lower bound by the properties of $g_T$ and $h_T$ given in 
Lemma \ref{lem:positive-test-family},
\[
p^{2Br-\epsilon/2} \cdot  N_{\mathrm{Dir}}(r;p)  \, \ll_{B,\epsilon,r}\, \sum_{\substack{\chi\bmod p\\
                \chi(-1)=1,\ \chi\neq \mathbf 1}}
\ \sum_j h_T(r_{\chi,j}).
\]
Combining these two inequalities, \eqref{eq:densitycount} follows easily.
\end{proof}

\begin{remark}
It follows from \cite[Theorems~1.5 and~1.9]{Humphries} that
\[
 N_{\mathrm{Dir}}(r;p)\ll_{\varepsilon,r}p^{2-8r+\varepsilon}.
\]
Under the conjectural values \(C \leq 2\) and \(\nu=0\) in \eqref{eq:arithmetic-input}, Proposition~\ref{prop:conditional-dirichlet-gap} gives
\[
 N_{\mathrm{Dir}}(r;p)\ll_{B,\varepsilon,r}p^{2-2Br+\varepsilon}.
\]
Thus \(B=4\) recovers Humphries' exponent in \cite{Humphries}, while every \(B>4\) gives a stronger density estimate in the nontrivial even Dirichlet sector. If \(B>1/r\), the exponent is negative and hence \(N_{\mathrm{Dir}}(r;p)=0\) for all sufficiently large \(p\). 
In particular, if \(B>64/7\), which is the Kim--Sarnak threshold of \cite{KimSarnak}, then \(N_{\mathrm{Dir}}(r;p)=0\) for all sufficiently large \(p\), for every fixed \(r\) satisfying \(1/B<r<7/64\).
\end{remark}

\subsection{Proof of Proposition \ref{hyp:trace-remainder}} 
\label{subsec:Dirichlet-remainder}

We now bound the remainder \eqref{eq:R-dir} using the vector-valued formulation of the trace formula. 
We use the weight-zero Selberg trace formula with the
finite-dimensional unitary representation
$\rho_{\mathrm{Dir},\Gamma}$; see
\cite[Chap.~VIII, \S4, Theorem~4.2, pp.~314--315]{HejhalII}.
Its geometric side consists of the identity, hyperbolic,
elliptic, and parabolic contributions, which we denote by
\[
I_{\mathrm{Dir}}(T), \quad H_{\mathrm{Dir}}(T),\quad
E_{\mathrm{Dir}}(T), \quad P_{\mathrm{Dir}}(T),
\]
respectively. For the vector-valued identity and hyperbolic
orbital integrals, see also
\cite[Chap.~III, \S4, p.~345~ff.]{HejhalI};
the elliptic contribution is treated in
\cite[Chap.~III, \S5, p.~350~ff.]{HejhalI}.
The parabolic calculation, including the gamma-factor
integral and the logarithms associated with nontrivial cusp
eigenvalues, appears in
\cite[Chap.~VIII, \S4, Eq.~(4.13), p.~313]{HejhalII}.

The continuous-spectrum contribution is expressed through
the logarithmic derivative of the determinant of the scattering matrix
$\Phi_{\chi}(s)$ in the $\chi$-sector, using
\[
 \frac{d}{ds}\log\det\Phi_\chi(s)
 =
 \operatorname{tr}\bigl(\Phi_\chi'(s)\Phi_\chi(s)^{-1}\bigr).
\]
The scattering matrix is defined in
\cite[Chap.~VIII, \S3, Definition~3.8, p.~281]{HejhalII}.
We retain the notation $C_{\mathrm{Dir}}(T)$ for
the scattering integral and include the terms proportional
to $h_T(0)$ in $P_{\mathrm{Dir}}(T)$.
After separating the primitive hyperbolic contribution,
the remaining terms, with the signs prescribed by the trace
formula, constitute $R_{\mathrm{Dir}}(T)$.

For each nontrivial even character $\chi$ modulo $p$, the representation $\operatorname{Ind}_{\Gamma_0(p)}^{\Gamma(1)}\chi_{\Gamma_0(p)}$ in \eqref{eq:induced-rep-Gamma}
has dimension $[\Gamma(1):\Gamma_0(p)]=p+1$. Summing over the $(p-3)/2$ such characters gives
\begin{equation}
 D_p:=\dim\rho_{\mathrm{Dir},\Gamma}
 =\frac{p-3}{2}(p+1)\ll p^2.
 \label{eq:dimension-rho-Dir}
\end{equation}
In particular, by unitarity,
\begin{equation}
 \bigl|\operatorname{tr}\rho_{\mathrm{Dir},\Gamma}(\gamma)\bigr|
 \le D_p\ll p^2
 \qquad(\gamma\in\Gamma(1)).
 \label{eq:trace-rho-Dir-trivial}
\end{equation}

We next estimate the contribution of the scattering part. The following lemma gives the bound needed for our family of test functions.

\begin{lemma}
\label{lem:Dirichlet-scattering}
For the following scattering contribution
\begin{equation*}
 C_{\mathrm{Dir}}(T)
 :=-\frac{1}{4\pi}
 \sum_{\substack{\chi\ ({\rm mod}\ p)\\
          \chi(-1)=1,\ \chi\ne\mathbf 1}}
 \int_{\mathbb R}h_T(r)\,
 \operatorname{tr}\!\left(
   \Phi_\chi'\!\left(\tfrac12+ir\right)
   \Phi_\chi\!\left(\tfrac12+ir\right)^{-1}
 \right)\,dr,
 \label{eq:continuous-term-definition}
\end{equation*}
where 
$h_T$ denotes the function constructed in Lemma~\ref{lem:positive-test-family}, we have
\begin{equation}\label{eq:continuous-coarse}
 C_{\mathrm{Dir}}(T)
 \ll_{\kappa}p^2\, T \log(2p).
\end{equation}
\end{lemma}

\begin{proof}
For each nontrivial even character $\chi$ modulo $p$, let
$\rho_{\chi,\Gamma}:=\rho_\chi\circ\pi_p$.
We use the scattering matrix in the induced cusp normalization
associated with $\rho_{\chi,\Gamma}$ in the vector-valued trace formula.
For a nonprincipal even character $\psi$ modulo $p$, put
\[
\Lambda(s,\psi)
:=
\left(\frac p\pi\right)^{s/2}
\Gamma\!\left(\frac s2\right)L(s,\psi),
\qquad
R_\psi(s):=\frac{\Lambda(2s-1,\psi)}{\Lambda(2s,\psi)}.
\]
The constant-term and functional-equation formulas for the
Eisenstein series give
\[
\det\Phi_\chi(s)=c_\chi R_\chi(s)R_{\bar\chi}(s),
\]
where $c_\chi\ne0$ is independent of $s$; see
\cite[Propositions~4.1--4.2 and equation~(6.2)]{Young}.
Here the passage from the width-one scalar cusp basis to the induced
cusp basis removes the factor $p^{1-2s}$ in the scalar scattering
determinant.

Logarithmic differentiation gives
\[
\operatorname{tr}\bigl(\Phi_\chi'(s)\Phi_\chi(s)^{-1}\bigr)
=
2\sum_{\psi\in\{\chi,\bar\chi\}}
\left(
\frac{\Lambda'}{\Lambda}(2s-1,\psi)
-
\frac{\Lambda'}{\Lambda}(2s,\psi)
\right),
\]
where the sum is taken with multiplicity, also when $\chi=\bar\chi$.
The functional equation and conjugation identity imply
\[
\frac{\Lambda'}{\Lambda}(2ir,\psi)
=
-\overline{
\frac{\Lambda'}{\Lambda}(1+2ir,\psi)
}.
\]
Consequently, on the unitary axis,
\[
\operatorname{tr}\bigl(
\Phi_\chi'(\tfrac12+ir)\Phi_\chi(\tfrac12+ir)^{-1}
\bigr)
=
-4\sum_{\psi\in\{\chi,\bar\chi\}}
\operatorname{Re}\frac{\Lambda'}{\Lambda}(1+2ir,\psi).
\]
Thus it suffices to estimate the integrals against these real
logarithmic derivatives.
The construction in Lemma~\ref{lem:positive-test-family} gives, for every $N\geq 0$,
\begin{equation}
 |h_T(r)|\ll_{N,\kappa}T(1+T|r|)^{-N}.
 \label{eq:hT-Schwartz-bound}
\end{equation}

Put
\(F_\psi(r):=\operatorname{Re}(\Lambda'/\Lambda)(1+2ir,\psi)\).
For primitive nonprincipal \(\psi\), write its nontrivial zeros as
\(\rho=\beta+i\gamma\), counted with multiplicity. The standard local
zero-counting estimate gives
\[
 n_\psi(u):=\#\{\rho:u<\gamma\le u+1\}
 \ll \log\bigl(\, p\,(2+|u|)\, \bigr),
 \qquad u\in\mathbb R.
\]
By the Hadamard product,
\[
 F_\psi(r)=
 \sum_\rho
 \frac{1-\beta}{(1-\beta)^2+(2r-\gamma)^2}\ge0.
\]
Each summand has integral \(\pi/2\) over \(\mathbb R\). Applying the
local zero count to nearby zeros and summing the quadratically
decaying tails therefore yields
\(\int_n^{n+1}F_\psi(r)\,dr
 \ll\log(p(2+|n|))\).
Consequently, \eqref{eq:hT-Schwartz-bound}, with \(N=3\), implies
\[
\begin{aligned}
 \left|\int_{\mathbb R}h_T(r)F_\psi(r)\,dr\right|
 &\le \sum_{n\in\mathbb Z}
 \sup_{r\in[n,n+1]}|h_T(r)|
 \int_n^{n+1}F_\psi(r)\,dr\\
 &\ll_\kappa T\sum_{n\in\mathbb Z}
 \frac{\log(p(2+|n|))}{(1+|n|)^3}
 \ll_\kappa T\log(2p).
\end{aligned}
\]
There are \((p-3)/2\) nonprincipal even characters modulo \(p\). For each \(\chi\), the factorization of \(\det\Phi_\chi(s)\) gives two logarithmic-derivative contributions, corresponding to \(R_\chi\) and \(R_{\bar\chi}\).
The total number of factor occurrences is therefore
$p-3$. Summing the preceding integral bound gives
\[
|C_{\mathrm{Dir}}(T)|
\ll_\kappa (p-3)T\log(2p)
\ll_\kappa p^2T\log(2p),
\]
which proves \eqref{eq:continuous-coarse}.
\end{proof}

We restate Proposition~\ref{hyp:trace-remainder} here
for convenience.

\medskip
\noindent\textbf{Proposition~\ref{hyp:trace-remainder} (restated).}
For every fixed $B>0$ and every $\epsilon>0$, one has
\begin{equation}
\label{eq:Dirichlet-remainder-bound}
|R_{\mathrm{Dir}}(T)| \ll_{B,\epsilon,\kappa} p^{2+\epsilon}
\qquad \text{for} \quad 1\leq T\leq 2B\log p.
\end{equation}

\begin{proof}
The sum of the scalar trace formulas over the nontrivial even characters is
the vector-valued Selberg trace formula for $\rho_{\mathrm{Dir}}$.  We bound
in turn the identity, elliptic, nonprimitive hyperbolic, parabolic, and
continuous terms.  The primitive hyperbolic term with exponent $m=1$ is
not included here, since it is precisely $H_{\mathrm{Dir},1}(T)$.

First, the identity contribution is
\[
 \frac{D_p\operatorname{vol}(\Gamma(1)\backslash\mathbb H)}{4\pi}
 \int_{\mathbb R}r\tanh(\pi r)h_T(r)\,dr.
\]
By \eqref{eq:hT-moments-remainder} and \eqref{eq:dimension-rho-Dir}, its absolute value is $O_\kappa(p^2)$.

Put \(\rho=\rho_{\mathrm{Dir},\Gamma}\). Since
\(\rho(-I)=I_{D_p}\), we compute the elliptic contribution on
\(\Gamma(1)/\{\pm I\}\). Let
\[
 S=\begin{pmatrix}0&-1\\1&0\end{pmatrix},
 \qquad
 R=\begin{pmatrix}0&-1\\1&1\end{pmatrix}.
\]
Their images in \(\Gamma(1)/\{\pm I\}\) have orders \(2\) and \(3\), respectively,
and \(S,R,R^2\) represent the three nontrivial elliptic conjugacy
classes in the projective quotient. Since \(h_T\) is even, their
contributions are
\[
\begin{aligned}
 E_{\mathrm{Dir}}(T)
 ={}&\frac{\operatorname{tr}\rho(S)}{8}
       \int_{\mathbb R}\frac{h_T(r)}{\cosh(\pi r)}\,dr\\
 &+\frac{\operatorname{tr}\rho(R)+\operatorname{tr}\rho(R^2)}
          {6\sqrt{3}}
       \int_{\mathbb R}
       h_T(r)\frac{\cosh(\pi r/3)}{\cosh(\pi r)}\,dr.
\end{aligned}
\]
The coefficients $1/8$ and $1/(6\sqrt{3})$ arise from
the general elliptic factor
$(4n\sin(\pi k/n))^{-1}$ for the $k$-th power
of an elliptic element of projective order $n$.
The two kernels satisfy
\[
 \frac{1}{\cosh(\pi r)}\le 2e^{-\pi|r|},
 \qquad
 \frac{\cosh(\pi r/3)}{\cosh(\pi r)}
 \le 2e^{-2\pi|r|/3},
\]
and both are bounded by \(1\). Thus
\eqref{eq:trace-rho-Dir-trivial} gives
\[
 |E_{\mathrm{Dir}}(T)|
 \le D_p\left(\frac18+\frac{1}{3\sqrt3}\right)
       \int_{\mathbb R}|h_T(r)|\,dr.
\]
Combining this with the identity contribution above, using
\(|r\tanh(\pi r)|\le |r|\), and applying
\eqref{eq:hT-moments-remainder}, we obtain
\begin{equation}
 |I_{\mathrm{Dir}}(T)|+|E_{\mathrm{Dir}}(T)|
 \ll D_p\int_{\mathbb R}(1+|r|)|h_T(r)|\,dr
 \ll_\kappa D_p
 \ll_\kappa p^2.
 \label{eq:identity-elliptic-bound}
\end{equation}

We next consider the nonprimitive hyperbolic terms.
Every nonprimitive hyperbolic element can be written as
$\gamma^m$, where $\gamma$ is primitive and $m\geq 2$.
We denote their total contribution by
$H_{\mathrm{Dir},\geq 2}(T)$.
Using
\eqref{eq:trace-rho-Dir-trivial}, their total absolute value is at most a
constant times
\begin{equation*}
 D_p\sum_{\{\gamma\}}\sum_{\substack{m\geq2\\ m \, \ell(\gamma)\leq T}}
 \frac{\ell(\gamma)}{2\sinh(m \, \ell(\gamma)/2)}
 |g_T(m \, \ell(\gamma))|.
 \label{eq:imprimitive-start}
\end{equation*}
Let
\[
\ell_0=2\log\frac{3+\sqrt5}{2}
\]
be the minimum length of a primitive closed geodesic on
$\Gamma(1)\backslash\mathbb H$. The prime geodesic theorem gives
\[
\#\{\gamma\in\mathcal P(\Gamma(1)):
       \ell_0+j\le\ell(\gamma)<\ell_0+j+1\}
\ll \frac{e^j}{1+j},
\qquad j\ge0.
\]
Since $\ell(\gamma)\ge\ell_0$ and $m\ge2$, we have, uniformly in $m$,
\[
\frac{\ell(\gamma)}{2\sinh(m\ell(\gamma)/2)}
=
\frac{\ell(\gamma)e^{-m\ell(\gamma)/2}}
     {1-e^{-m\ell(\gamma)}}
\ll \ell(\gamma)e^{-m\ell(\gamma)/2}.
\]
Together with $\|g_T\|_\infty\ll_\kappa1$, this shows that the
contribution $H_{\mathrm{Dir},m}(T)$ of the $m$-th powers satisfies
\[
|H_{\mathrm{Dir},m}(T)|
\ll_\kappa
D_p e^{-m\ell_0/2}
\sum_{\substack{j\ge0\\ \ell_0+j\le T/m}}
e^{-(m/2-1)j}.
\]
For $m=2$, this is $O_\kappa(D_p(1+T))$. For $m\ge3$, the
geometric sum is bounded uniformly by $(1-e^{-1/2})^{-1}$, and hence
\[
\sum_{m\ge3}|H_{\mathrm{Dir},m}(T)|
\ll_\kappa D_p\sum_{m\ge3}e^{-m\ell_0/2}
\ll_\kappa D_p.
\]
Consequently,
\begin{equation}\label{eq:imprimitive-final}
|H_{\mathrm{Dir},\ge2}(T)|
\ll_\kappa D_p(1+T)
\ll_\kappa p^2(1+T).
\end{equation}

We next estimate the parabolic contribution. Let
$\tau=\left(\begin{smallmatrix}1&1\\0&1\end{smallmatrix}\right)$
generate the stabilizer of the unique cusp, modulo the center.
The scattering matrix $\Phi_\chi(s)$ acts on the
$\rho_{\chi,\Gamma}(\tau)$-invariant subspace
\[
\ker(\rho_{\chi,\Gamma}(\tau)-I),
\]
which has dimension $d_\chi^{\mathrm{Eis}}$.
Since $\rho$ is the direct sum of these representations,
\[
M_p:=
\sum_{\substack{\chi\bmod p\\\chi(-1)=1,\ \chi\ne\mathbf 1}}
d_\chi^{\mathrm{Eis}}
=\dim\ker(\rho(\tau)-I)\le D_p.
\]
With the scattering integral treated separately, the
parabolic contribution consists of fixed numerical multiples
of
\[
 M_p g_T(0),\qquad
 M_p\int_{\mathbb R}h_T(r)
       \frac{\Gamma'}{\Gamma}(1+ir)\,dr,\qquad
 g_T(0)\sum_{\lambda\ne1}m_\lambda\log|1-\lambda|,
\]
together with the terms involving
\[
 h_T(0)M_p
 \quad\text{and}\quad
 h_T(0)
 \sum_{\substack{\chi\bmod p\\\chi(-1)=1,\ \chi\ne \mathbf 1}}
 \operatorname{tr}\Phi_\chi(\tfrac12).
\]
Here $\lambda$ in the third sum runs over the nontrivial eigenvalues of
$\rho(\tau)$, with multiplicity $m_\lambda$, so that
$\sum_{\lambda\ne1}m_\lambda=D_p-M_p$. 
Since $\rho(\tau)^p=I$, every such $\lambda$ is a nontrivial
$p$-th root of unity. Therefore
\[
 2\sin(\pi/p)\le |1-\lambda|\le2,
 \qquad
 \bigl|\log|1-\lambda|\bigr|\ll\log(2p),
\]
and the logarithmic term is
$O_\kappa((D_p-M_p)\log(2p))$.
Moreover,
$\bigl|\Gamma'/\Gamma(1+ir)\bigr|\ll\log(2+|r|)$,
so \eqref{eq:hT-moments-remainder} gives
\[
 \left|\int_{\mathbb R}h_T(r)
       \frac{\Gamma'}{\Gamma}(1+ir)\,dr\right|
 \ll_\kappa1.
\]
Now $\Phi_\chi(\tfrac12)^2=I$ implies
$|\operatorname{tr}\Phi_\chi(\tfrac12)|
\le d_\chi^{\mathrm{Eis}}$. Using
$g_T(0)\ll_\kappa1$ and $h_T(0)\ll_\kappa T$, we conclude
\begin{equation} \label{eq:parabolic-final}
 |P_{\mathrm{Dir}}(T)|
 \ll_\kappa
 M_p(1+T)+(D_p-M_p)\log(2p)
 \ll_\kappa p^2\bigl(T+\log(2p)\bigr).
\end{equation}

Combining
\eqref{eq:identity-elliptic-bound},
\eqref{eq:imprimitive-final},
\eqref{eq:parabolic-final}, and
\eqref{eq:continuous-coarse}, we obtain
\begin{equation}
 |R_{\mathrm{Dir}}(T)|
 \ll_{\kappa}p^2(1+T)\,\log(2p).
 \label{eq:remainder-before-log-range}
\end{equation}
Finally, if $1\leq T\leq2B\log p$, then
\[
 (1+T)\, \log(2p)
 \ll_{B}(\log(2p))^{2}
 \ll_{B,\epsilon}p^\epsilon.
\]
Substitution into \eqref{eq:remainder-before-log-range} proves
\eqref{eq:Dirichlet-remainder-bound}.
\end{proof}

\section*{Acknowledgments} 
%We would like to thank -----------------------------------------------------------------.
C.-h. Lee was supported by a KIAS Individual Grant (AP067301) via the Center for Artificial Intelligence and Natural Sciences at Korea Institute for Advanced Study.
K.-H. Lee was supported by NSF grant DMS-2617278. J. Park was supported by a KIAS Individual Grant (MG022804) via the School of Mathematics at Korea Institute for Advanced Study.

\section*{Code Availability Statement}
%https://www.sciencedirect.com/journal/journal-of-number-theory/publish/guide-for-authors#about-the-journal-article-types
The code and data used to generate the figures are available at
\url{https://github.com/chlee-0/character-twisted-class-number-sums}.

\section*{Declaration of generative AI use}

During the preparation of this manuscript, the authors used OpenAI ChatGPT to explore mathematical ideas, to develop and debug computational code, and to assist with editorial and formatting checks. 
The authors reviewed the AI-generated suggestions and computational results.

\bigskip

\noindent
\textsc{Chul-hee Lee}\\
Center for Artificial Intelligence and Natural Sciences,\\
Korea Institute for Advanced Study,\\
85 Hoegiro, Dongdaemun-gu, Seoul 02455, Republic of Korea\\
\textit{Email:} \texttt{chlee@kias.re.kr}

\medskip

\noindent
\textsc{Kyu-Hwan Lee}\\
Department of Mathematics, University of Connecticut,\\
341 Mansfield Road, U1009,\\
Storrs, Connecticut 06269-1009, USA
\smallskip
\\
Center for Artificial Intelligence and Natural Sciences,\\
Korea Institute for Advanced Study,\\
85 Hoegiro, Dongdaemun-gu, Seoul 02455, Republic of Korea\\\textit{Email:} \texttt{khlee@math.uconn.edu}

\medskip

\noindent
\textsc{Jinsung Park}\\
School of Mathematics, Korea Institute for Advanced Study,\\
85 Hoegiro, Dongdaemun-gu, Seoul 02455, Republic of Korea\\
\textit{Email:} \texttt{jinsung@kias.re.kr}

\end{document}